\documentclass[twoside, 11pt]{article}

\usepackage{jmlr2e}
\usepackage[linesnumbered]{algorithm2e}
\SetKw{Continue}{continue}
\jmlrheading{18}{2017}{1-31}{2/17; Revised 7/17}{12/17}{17-101}{K. S. Sesh Kumar and Francis Bach}

\ShortHeadings{Active-set Methods for  Submodular Minimization Problems}{Kumar and Bach}
\firstpageno{1}

\usepackage{times}
\usepackage{rotating}
\usepackage{epsfig}
\usepackage{graphicx}
\usepackage{amsmath}
\usepackage{amsthm}
\usepackage{amssymb}
\usepackage{multirow}
\usepackage{tabularx}
\usepackage{booktabs}
\usepackage{url,color}

\newcommand{\BEAS}{\begin{eqnarray*}}
\newcommand{\EEAS}{\end{eqnarray*}}
\newcommand{\BEA}{\begin{eqnarray}}
\newcommand{\EEA}{\end{eqnarray}}
\newcommand{\BEQ}{\begin{equation}}
\newcommand{\EEQ}{\end{equation}}
\newcommand{\BIT}{\begin{itemize}}
\newcommand{\EIT}{\end{itemize}}
\newcommand{\BNUM}{\begin{enumerate}}
\newcommand{\ENUM}{\end{enumerate}}
\newcommand{\BA}{\begin{array}}
\newcommand{\EA}{\end{array}}

\newcommand{\argmin}{\mathop{\rm argmin}}
\newcommand{\argmax}{\mathop{\rm argmax}}

\newcommand{\rb}{\mathbb{R}}
\newcommand{\lova}{Lov\'asz }

\renewcommand{\labelitemi}{$-$}

\newcommand{\mysec}[1]{Section~\ref{sec:#1}}
\newcommand{\myapp}[1]{Appendix~\ref{app:#1}}
\newcommand{\eq}[1]{Eq.~\ref{eq:#1}}
\newcommand{\myfig}[1]{Figure~\ref{fig:#1}}
\makeatletter
\newcommand{\subalign}[1]{%
  \vcenter{%
    \Let@ \restore@math@cr \default@tag
    \baselineskip\fontdimen10 \scriptfont\tw@
    \advance\baselineskip\fontdimen12 \scriptfont\tw@
    \lineskip\thr@@\fontdimen8 \scriptfont\thr@@
    \lineskiplimit\lineskip
    \ialign{\hfil$\m@th\scriptstyle##$&$\m@th\scriptstyle{}##$\crcr
      #1\crcr
    }%
  }
}
\makeatother

\allowdisplaybreaks

\begin{document}
\title{Active-set Methods for  Submodular Minimization Problems}
\author{\name K. S. Sesh Kumar \email s.karri@imperial.ac.uk \\
        \addr Imperial College London, Department of Computing\\
	180 Queen's Gate, London SW7 2AZ
        \AND 
        \name Francis Bach \email francis.bach@ens.fr\\
	\addr INRIA - Sierra Project-team \\ 
	D{\'e}partement d{'}Informatique de l{'}Ecole Normale Sup{\'e}rieure (UMR CNRS/ENS/INRIA)\\
        2, rue Simone Iff \\
	75012 Paris, France}

\editor{Sebastian Nowozin}

\maketitle

\begin{abstract}%
We consider the submodular function minimization~(SFM) and the quadratic minimization problems regularized by the \lova extension of the submodular function. These optimization problems are intimately related; for example,  min-cut problems and total variation denoising problems, where the cut function is submodular and its \lova extension is given by the associated total variation. When a quadratic loss is regularized by the total variation of a cut function, it thus becomes a total variation denoising problem and we use the same terminology in this paper for ``general'' submodular functions.
 We propose a new active-set algorithm for total variation denoising with the assumption of an oracle that solves the corresponding SFM problem. This can be seen as local descent algorithm over ordered partitions with explicit convergence guarantees. It is more flexible than the existing algorithms with the ability for warm-restarts using the solution of a closely related problem. Further, we also consider the case when a submodular function can be decomposed into the sum of two submodular functions $F_1$ and $F_2$ and assume SFM oracles for these two functions. We propose a new active-set algorithm for total variation denoising (and hence SFM by thresholding the solution at zero). This algorithm also performs local descent over ordered partitions and its ability to warm start considerably improves the performance of the algorithm. In the experiments, we compare the performance of the proposed algorithms with  state-of-the-art algorithms, showing that it reduces the calls to SFM oracles.
\end{abstract}


\section{Introduction}

 Submodular optimization problems such as total variation denoising and submodular function minimization are convex optimization problems which are common in computer vision, signal processing and machine learning~\citep{fot_submod}, with notable applications to graph cut-based image segmentation~\citep{boykov2001fast}, sensor placement~\citep{krause11submodularity}, or document summarization~\citep{lin2011-class-submod-sum}. 
 
Let $F$ be a normalized submodular function defined on $V = \{1, \dots, n\}$, i.e., $F:2^V \to \rb$ such that $F(\varnothing) = 0$ and an $n$-dimensional real vector $u$, i.e., $u \in \rb^n$. In this paper, we consider the submodular function minimization (SFM) problem,
\BEQ
\min_{A \subset V} F(A) - u(A),
\label{eq:sfm}
\EEQ
where we use the convention $u(A) = u^\top 1_A$ and $1_A \in \{0,1\}^n$ is the indicator vector of the set $A$. Note that general submodular functions can always be decomposed into a normalized submodular function, $F$, i.e., $F(\varnothing) = 0$ and a modular function  $u$ \citep[see][]{fot_submod}.

Let $f$ be the \lova extension of the submodular function $F$. Let us consider the following continuous optimization problem
\BEQ
\min_{w \in [0, 1]^n}f(w) - u^\top w.
\label{eq:relaxsfm}
\EEQ
As a consequence of submodularity, the discrete and continuous optimization problems in \eq{sfm} and \eq{relaxsfm} respectively have the same optimal solutions~\citep{lovasz1982submodular}. Let us consider another related continuous optimization problem 
\BEQ
\min_{w \in \rb^n} f(w) - u^\top w + \textstyle \frac{1}{2} \| w  \|_2^2.
\label{eq:tv}
\EEQ
If $F$ is a cut function in a weighted undirected graph, then $f$ is its associated total variation, hence the denomination of  \emph{total variation denoising} (TV) problem for \eq{tv}, which we use in this paper---since it is equivalent to minimizing $\textstyle \frac{1}{2} \| u - w\|_2^2 + f(w)$. The unique solution of the total variation denoising in \eq{tv} can be used to obtain the solution of the SFM problem in \eq{sfm} by thresholding at $0$. Conversely, we may obtain the optimal solution of the total variation denoising in \eq{tv} by solving a series of SFM problems using divide-and-conquer strategy.

{\it Relationship with existing work.} Generic algorithms to optimize SFM in \eq{sfm} or TV in \eq{tv} problems which only access $F$ through function values, e.g., subgradient descent or min-norm-point algorithm~\citep{fujishige1984submodular}, are too slow without any assumptions~\citep{fot_submod}, as  for signal processing applications, high precision is typically required (and often the exact solution).

 For decomposable problems, i.e., when $F = F_1 + \cdots + F_r$, where each $F_j$ is ``simple'', some algorithms use more powerful oracles than function evaluations, improving the running times. These powerful oracles include SFM oracles that can solve the SFM problem of simple submodular function, $F_j$ given by
\BEQ
\min_{A \subset V} F_j(A) - u_j(A),
\label{eq:SFMj}
\EEQ
where $u_j \in \rb^n$. The other set of powerful oracles are total variation or TV oracles, that solve TV problems of the form
\BEQ
\min_{w \in \rb^n} f_j(w) - u_j^{\top}w  +  \textstyle \frac{1}{2} \| w  \|_2^2,
\label{eq:TVj}
\EEQ
where $u_j \in \rb^n$. Note that, in general, the exact total variation oracles are at most $O(n)$ times more expensive than their respective SFM oracle as they solve all SFM problems
\BEQ
\min_{A \subset V} F_j(A) - u_j(A) + \lambda |A|,
\EEQ
for all $\lambda \in \rb$, which have at most $O(n)$ unique solutions. For more details refer to ~\citet{fujishige80} and \citet{fot_submod}. Here, $|A|$ denotes the cardinality of the set $A$. There does exist a subclass of submodular functions (cut functions and other submodular functions that can be written in form of cuts) whose total variation oracles are only $O(1)$ times more expensive than the corresponding SFM oracles but are still too expensive in practice.

 \citet{stobbe11} used SFM oracles instead of function value oracles but their algorithm remains slow in practice. However, when total variation oracles for each $F_j$ are used, they become competitive~\citep{komodakis2011mrf,seshTV,treesubmod}. Therefore, our goal is to design fast optimization strategies using only efficient SFM oracles for each function $F_j$ rather than their expensive TV oracles~\citep{seshTV,treesubmod} to solve the SFM and TV denoising problems of $F$ given by \eq{sfm} and \eq{tv} respectively. An algorithm was proposed by ~\citet{loic2016} to search over partition space for solving \eq{tv} with the unary terms $(- u^{\top}w)$ replaced by a convex differentiable function but it applies only to functions $F$, which are cut functions. 

 In this paper, we exploit the polytope structure of these non-smooth optimization problems with exponentially many constraints, i.e., $2^n$, where each face of the constraint set is indexed by an ordered partition of the underlying ground set $V = \{1,\dots,n\}$. The main insight of this paper is that given the main polytope associated with a submodular function (namely the base polytope described in \mysec{reviewsubmod}) and an ordered partition, we may uniquely define a tangent cone of the polytope. Further, orthogonal projections onto the tangent cone may be done efficiently by isotonic regressions~\citep{best1990active}. The time needed is linear in the number of elements of the ordered partition used to define the tangent cone. We need SFM oracles only to check the optimality of the ordered partition. Given the orthogonal projection $s$ onto the tangent cone, if the minimum of $F(A) - s(A)$ with respect to $A \subseteq V$ is positive then it is optimal. If it is not optimal, it gives us the violating constraints in the form of active-sets that enable us to generate a new ordered partition among the exponentially many ordered partitions.
 
{\it Contributions.}
We make two main contributions: 
\begin{list}{\labelitemi}{\leftmargin=1.7em}
\addtolength{\itemsep}{-.215\baselineskip}
\item Given a submodular function $F$ with an SFM oracle,  we propose a new active-set algorithm for total variation denoising in \mysec{nondecompalg}, which is more efficient and flexible than existing ones. This algorithm may be seen as a local descent algorithm over ordered partitions. It has the additional advantage of allowing warm restarts, which will be beneficial when we have to solve a large number of total variation denoising problems as  shown in \mysec{experiments}.

\item Given a decomposition of $F = F_1 + F_2$, with available SFM oracles for each $F_j$, we propose an active-set algorithm for total variation denoising in \mysec{decompalg} (and hence for SFM by thresholding the solution at zero). These algorithms optimize over ordered partitions (one per function~$F_j$). Following~\citet{treesubmod} and \citet{seshTV}, they are also naturally parallelizable. Given that only SFM oracles are needed, it is much more flexible than the algorithms requiring a TV oracle, and allow more applications as shown in \mysec{experiments}.
\end{list}

\begin{figure}
\begin{center}
\begin{tabular}{cc}
  \includegraphics[width=0.48\textwidth]{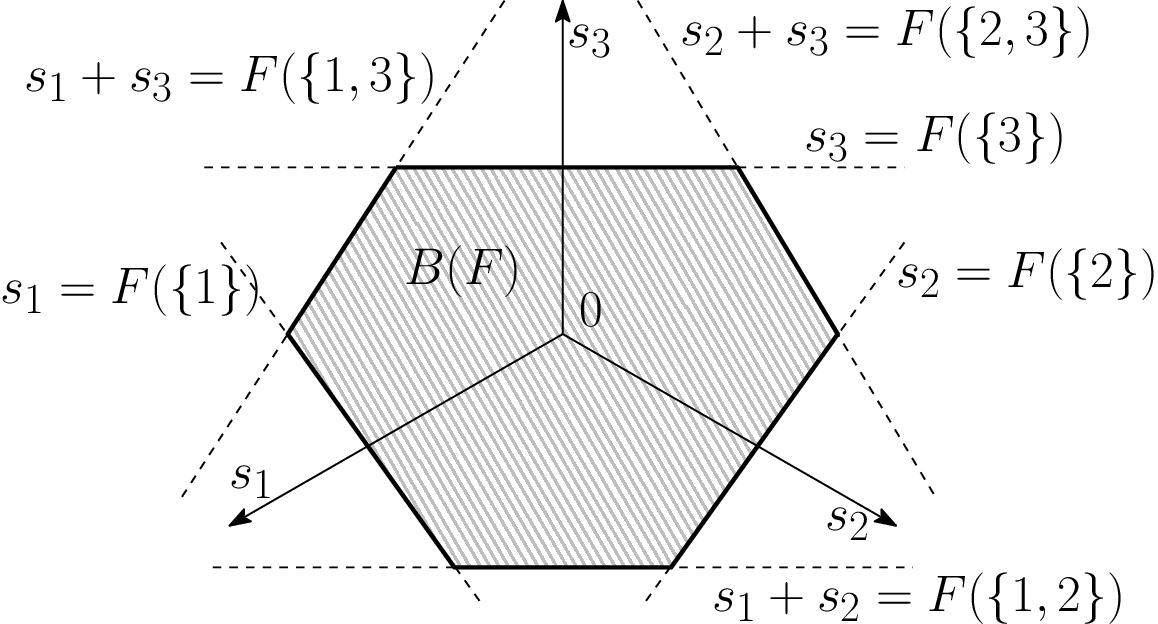} & 
  \includegraphics[width=0.48\textwidth]{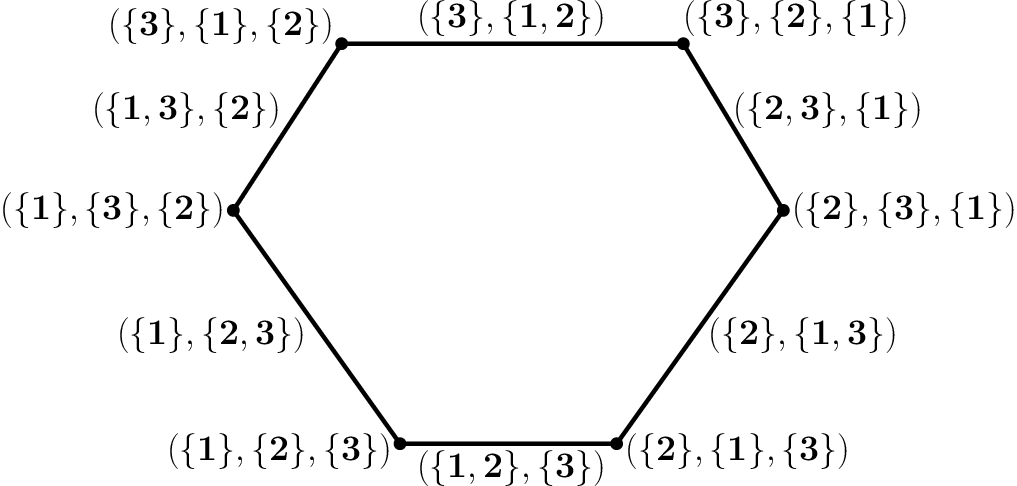}  \\
(a) & (b)
\end{tabular}
  \end{center}
  
  \vspace*{-.6cm}
  
  \caption{Base polytope for $n\!=\!3$. (a) definition from its supporting hyperplanes $\{s(A) = F(A)\}$. (b) each face (point or segment) of $B(F)$ is associated with an ordered partition. }
  \label{fig:base}
  \end{figure}

\section{Review of Submodular Analysis}
\label{sec:reviewsubmod}

A set-function $F: 2^V \to \rb$ is submodular if $F(A) + F(B) \geqslant F(A \cup B) + F (A \cap B)$ for any subsets $A,B$ of $V$. Our main motivating examples in this paper are cuts in a weighted undirected graph with weight function $a: V \times V \to \rb_+$, which can be defined as 
\BEQ
F(A) = \sum_{i < j} a(i, j)|1_{i \in A} - 1_{j \notin A}|,
\label{eq:cut}
\EEQ
where $i, j \in V$. Note that there are other submodular functions on which our algorithm works, e.g., concave function on the cardinality set, which \emph{cannot} 
be represented in the form of \eq{cut}~\citep{Ladicky2010,kolmogorov12}.
However, we use cut functions as a running example to better explain our algorithm as they are most widely studied and 
understood  among submodular functions.  We now  review the relevant concepts from submodular analysis \citep[for more details, see][]{fot_submod,fujishige2005submodular}.

{\it \lova extension and convexity.}
The power set $2^V$ is naturally identified with the vertices $\{0,1\}^n$ of the hypercube in $n$ dimensions (going from $A \subseteq V$ to $1_A \in \{0,1\}^n$). Thus, any set-function may be seen as a function $f$ on $\{0,1\}^n$. It turns out that $f$ may be extended to the full hypercube $[0,1]^n$ by piecewise-linear interpolation, and then to the whole vector space $\rb^n$. 

Given a vector $w \in \rb^n$, and given its \emph{unique} level-set representation as $w = \sum_{i=1}^m v_i 1_{A_i}$, with $(A_1,\dots,A_m)$ a partition of $V$ and $v_1 > \cdots > v_m$,  $f(w)$ is equal to $f(w) = \sum_{i=1}^m v_i \big[ F(B_i) - F(B_{i-1})\big]$, where $B_i = (A_1 \cup \cdots \cup A_i)$. For cut functions, the \lova extension happens to be equal to the \emph{total variation}, $f(w) = \sum_{i < j} a(i,j) |w_i-w_j|$, hence our denomination total variation denoising for the problem in \eq{tv}.

This extension is piecewise linear for any set-function~$F$. It turns out that it is convex if and only if $F$ is submodular~\citep{lovasz1982submodular}. Any piecewise linear convex function may be represented as the support function of a certain polytope $K$, i.e., as $f(w) = \max_{s \in K} w^\top s$~\citep{rockafellar97}. For the \lova extension of a submodular function, there is an explicit description of $K$, which we now review.

{\it Base polytope.} We define the \emph{base polytope} as
$$
B(F) = \\
\big\{ s \in \rb^n, \ s(V) = F(V), \ \forall A \subset V, s(A) \leqslant F(A) \big\}.
$$
Given that it is included in the affine hyperplane $\{s(V) = F(V)\}$, it is traditionally represented by the projection on that hyperplane (see  \myfig{base} (a)). A key result in submodular analysis is that the \lova extension is the support function of $B(F)$, that is, for any $w \in \rb^n$, 
\BEQ
f(w) = \sup_{s \in B(F)} w^\top s.
\label{eq:lova}
\EEQ
The maximizers  above may be computed in closed form from an ordered level-set representation of~$w$ using a greedy algorithm, which  (a) first sorts the elements of $w$ in decreasing order such that $w_{\sigma(1)} \geq \ldots \geq w_{\sigma(n)}$ where $\sigma$ represents the order of the elements in $V$; and (b) computes $s_{\sigma(k)} = F(\{ \sigma(1), \ldots, \sigma(k)\}) - F(\{ \sigma(1), \ldots, \sigma(k-1)\})$.

{\it SFM as a convex optimization problem.}
Another key result of submodular analysis is that minimizing a submodular function $F$ (i.e., minimizing the \lova extension~$f$ on $\{0,1\}^n$), is equivalent to minimizing the \lova extension $f$ on the full hypercube $[0,1]^n$ (a convex optimization problem). Moreover, with convex duality we have
\BEAS
\min_{A \subseteq V} F(A) - u(A) & = & \min_{w \in \{0,1\}^n} f(w) - u^{\top}w =  \min_{w \in [0,1]^n} f(w) - u^{\top}w \\
                                 & = & \min_{w \in [0,1]^n} \max_{s \in B(F)} s^{\top}w - u^{\top}w \\
                                 & = & \max_{s \in B(F)} \min_{w \in [0,1]^n} s^{\top}w - u^{\top}w
                                 = \max_{s \in B(F)} \sum_{i = 1}^n \min\{s_i - u_i, 0\}.
\EEAS
This dual problem allows to obtain certificates of optimality for the primal-dual pairs $w \in [0,1]^n$ and $s \in B(F)$ using the quantity
$$ \text{gap}(w,s) := f(w) - u^{\top}w - \sum_{i=1}^n \min\{ s_i - u_i, 0 \}, $$
which is always non-negative. It is equal to zero only at optimality and the corresponding $(w, s)$ form an optimal primal-dual pair.

{\it Total variation denoising as projection onto the base polytope.}
A consequence of the representation of $f$ as a support function leads to the following primal/dual pair~\citep[Sec.~8]{fot_submod}:
\BEA
\min_{w \in \rb^n} f(w)  -u^\top w + \textstyle \frac{1}{2} \| w  \|_2^2
& = & \min_{w \in \rb^n} \max_{s \in B(F)} s^\top w -u^\top w + \textstyle \frac{1}{2} \| w  \|_2^2 \text{ using \eq{lova}},  \nonumber \\
& = & \max_{s \in B(F)} \min_{w \in \rb^n} s^\top w -u^\top w + \textstyle \frac{1}{2} \| w  \|_2^2,  \nonumber \\
& = & \max_{s \in B(F)} -\textstyle \frac{1}{2}  \| s -u  \|_2^2, 
\label{eq:tvdual} 
\EEA
with $w = u-s$ at optimality.  Thus the TV problem is equivalent to the orthogonal projection of $u$ onto $B(F)$. 

{\it From TV denoising to SFM.}
The SFM problem in \eq{sfm} and the TV problem in \eq{tv} are tightly connected. Indeed, given the unique solution $w$ of the TV problem, we obtain a solution of $\min_{A \subseteq V} F(A) - u(A) $ by thresholding  $w$ at $0$, i.e., by taking $A = \{ i \in V, w_i \geqslant 0\}$~\citep{fujishige80}.

Conversely, one may solve the TV problem by an appropriate sequence of SFM problems. The original divide-and-conquer algorithm may involve $O(n)$ SFM problems~\citep{groenevelt1991two}. The extended algorithm of~\citet{treesubmod} can reach a precision $\varepsilon$ in $O(\log \textstyle \frac{1}{\varepsilon})$ iterations but can only get the exact solution in $O(n)$ oracles. Fast efficient algorithms are proposed to solve TV problems with $O(1)$ oracles~\citep{chambolle2009total,goldfarb2009} but are specific to cut functions on simple graphs (chains and trees) as they exploit the weight representation given by \eq{cut}. Our algorithm in \mysec{nondecompalg} is a generalization of the divide-and-conquer strategy for solving the TV problem with general submodular functions.

\section{Ordered Partitions and Isotonic Regression}

\label{sec:nondecompalg}

The main insight of this paper is (a) to consider the detailed face structure of the base polytope $B(F)$ and (b) to notice that for the outer approximation of $B(F)$ based on the tangent cone to a certain face, the orthogonal projection problem (which is equivalent to constrained TV denoising) may be solved efficiently using a simple algorithm, originally proposed to solve isotonic regression in linear time. This allows an explicit  efficient local search over ordered partitions.

\subsection{Outer Approximations of $B(F)$}

{\it Supporting hyperplanes.}
The base polytope is defined as the intersection of half-spaces $\{ s(A) \leqslant F(A)\}$, for $A \subseteq V$. Therefore, faces of $B(F)$ are indexed by subsets of the power set. As a consequence of submodularity~\citep{fot_submod,fujishige2005submodular}, the faces  of the base polytope $B(F)$ are characterized by ``ordered partitions'' $\mathcal{A} = (A_1,\dots,A_m)$ with  $V = A_1 \cup \cdots \cup A_m$. Then, a face of $B(F)$ is such that $s(B_i) = F(B_i)$ for all $B_i = A_1 \cup \cdots \cup A_i$, $i=1,\dots,m$. 
See \myfig{base}~(b) for the enumeration of faces for $n=3$ based on an enumeration of all ordered partitions. Such ordered partitions are associated to vectors $w = \sum_{i=1}^m v_i 1_{A_i}$ with $v_1 > \cdots > v_m$ with all solutions of $\max_{s \in B(F)} w^{\top}s$ being on the corresponding face.

\begin{figure}
\begin{center}
\begin{tabular}{cc}
\includegraphics[width=0.54\textwidth]{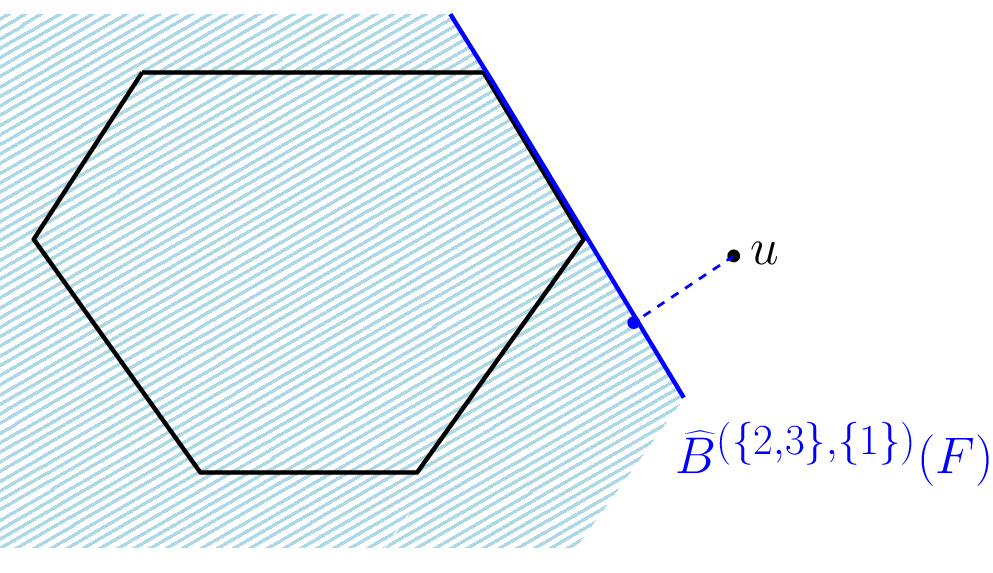} &
\includegraphics[width=0.46\textwidth]{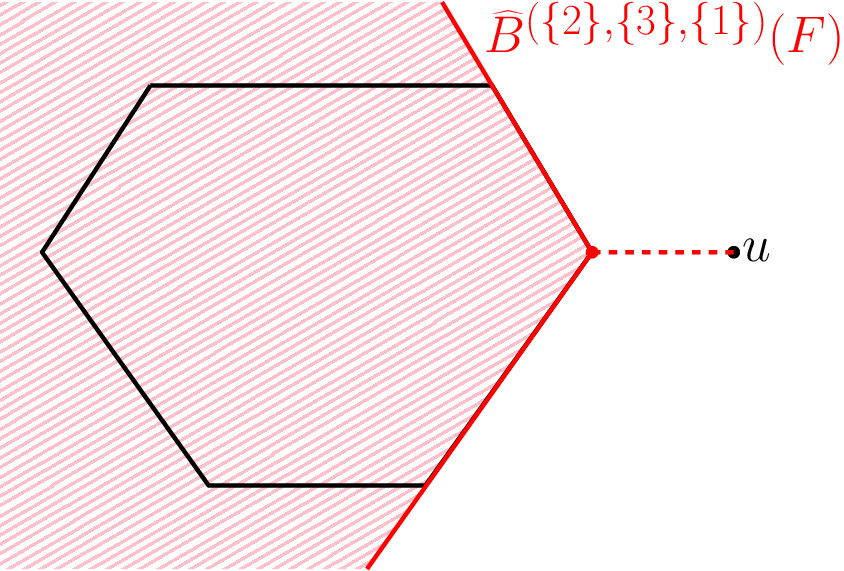}  
\end{tabular}
\end{center}

\vspace*{-.5cm}

\caption{Projection algorithm for a single polytope: first projecting on the outer approximation $\widehat{B}^{(\{2,3\},\{1\})}(F)$, with a projected element which is not in $B(F)$ (blue), then on $\widehat{B}^{(\{2\},\{3\},\{1\})}(F)$, with a projected element being the projection of $s$ onto $B(F)$ (red). }
\label{fig:outer}
\end{figure}

From a face of $B(F)$ defined by the ordered partition $\mathcal{A}$, we may define its \emph{tangent cone}  $\widehat{B}^\mathcal{A}(F)$ at this face as the set
\BEQ
\widehat{B}^{\mathcal{A}}(F) = \Big\{ s \in \rb^n,  s(V)=F(V), \forall i \in \{1,\dots,m-1\},  s(B_i) \leqslant F(B_i)\Big\}.  \label{eq:tangConeBase}
\EEQ
Since we have relaxed all the constraints unrelated to $\mathcal{A}$, these are outer approximations of $B(F)$ as illustrated in \myfig{outer} for two ordered partitions.

{\it Support function.} We may compute the support function of $\widehat{B}^\mathcal{A}(F)$, which is an upper bound on $f(w)$ since this set is an outer approximation of $B(F)$ as follows:
\vspace{-1em}
\BEAS
\sup_{ s\in \widehat{B}^{\mathcal{A}}(F)} w^\top s
&= &\sup_{s \in \rb^n} \inf_{ \lambda \in \rb^{m-1}_+ \times \rb} w^\top s - \sum_{i=1}^m \lambda_i ( s(B_i) - F(B_i) ) \mbox{ , using Lagrangian duality,} \\[-.1cm]
&= & \inf_{ \lambda \in \rb^{m-1}_+ \times \rb} \sup_{s \in \rb^n} s^\top \Big( w - \sum_{i=1}^m  ( \lambda_i + \cdots + \lambda_m ) 1_{A_i} \Big)\\[-.1cm]
&  &       \hspace*{4cm}             + \sum_{i=1}^m  ( \lambda_i + \cdots + \lambda_m ) \big[ F(B_i) - F(B_{i-1}) \big], \\
&= &\inf_{ \lambda \in \rb^{m-1}_+ \times \rb} \sum_{i=1}^m  ( \lambda_i + \cdots + \lambda_m ) \big[ F(B_i) - F(B_{i-1}) \big] \\[-.1cm]
&  &             \hspace*{4cm}             \mbox{ such that } w = \sum_{i=1}^m  ( \lambda_i + \cdots + \lambda_m ) 1_{A_i}.
\EEAS
Thus, by defining $v_i = \lambda_i + \cdots + \lambda_m$, which are decreasing, the support function is finite for $w$ having ordered level sets corresponding to the ordered partition $\mathcal{A}$ (we then say that $w$ is \emph{compatible} with $\mathcal{A}$). In other words, if $w = \sum_{i=1}^m v_i 1_{A_i}$, the support functions is equal to the \lova extension $f(w)$. Otherwise, when $w$ is not compatible with $\mathcal{A}$, the support function is infinite.

Let us now denote $\mathcal{W}^\mathcal{A}$ as a set of all weight vectors $w$ that are compatible with the ordered partition $\mathcal{A}$. This can be defined as

$$ \mathcal{W}^\mathcal{A} = \bigg\{w \in \rb^n \mid \exists v \in \rb^m, w = \sum_{i = 1}^m v_i 1_{A_i}, v_1 \geq \ldots \geq v_m \bigg\}.$$ 

Therefore, 
\BEQ
\sup_{ s\in \widehat{B}^{\mathcal{A}}(F)} w^\top s  = 
\begin{cases} 
f(w)  & \text{if } w \in \mathcal{W}^{\mathcal{A}} , \\
\infty & \text{otherwise.}
\end{cases} 
\label{eq:suppTangentCone}
\EEQ

\subsection{Isotonic Regression for Restricted Problems}
\label{sec:isotonic}

Given an ordered partition $\mathcal{A} = (A_1,\dots,A_m)$ of $V$, we consider the original TV problem restricted to $w$ in $\mathcal{W}^\mathcal{A}$.  
Since on this constraint set $f(w) = \sum_{i=1}^m v_i \big[ F(B_i) - F(B_{i-1})\big]$ is a linear function, this is equivalent to
\BEQ
\label{eq:isotonic}
\min_{v \in \rb^m}
\sum_{i=1}^m v_i \big[ F(B_i) - F(B_{i-1}) - u(A_i) \big] + \textstyle \frac{1}{2}\sum_{i=1}^m |A_i| v_i^2  \text{ such that }  v_1 \geqslant \cdots \geqslant v_m.
\EEQ
 This may be done by isotonic regression in complexity $O(m)$ by the weighted pool-adjacent-violator algorithm~\citep{best1990active}. Typically the solution~$v$ will have some values that are equal to each other, which corresponds to merging some sets $A_i$. If these merges are made, we now obtain a {\em basic ordered partition}\footnote{Given a submodular function $F$ and an ordered partition $\mathcal{A}$,  when the unique solution problem in \eq{isotonic} is such that $v_1 > \cdots > v_m$, we say that we $\mathcal{A}$ is a \emph{basic ordered partition} for $F-u$. Given any ordered partition, isotonic regression allows to compute a coarser partition (obtained by partially merging some sets) which is basic.} such that our optimal $w$ has  \emph{strictly decreasing} values. Because none of the constraints are tight, primal stationarity leads to explicit values of $v$ given by  $v_i =  {u(A_i)}/{|A_i|} -  {(F(B_i) - F(B_{i-1}))}/{|A_i|}$, i.e., given~$\mathcal{A}$, the  exact  solution of the TV problem may be obtained in closed form.

{\it Dual interpretation.}
\eq{isotonic} is a constrained TV denoising problem that minimizes the cost function in \eq{tv} but with the constraint that weights are compatible with 
the ordered partition $\mathcal{A}$, i.e.,  $\min_{w \in \mathcal{W}^\mathcal{A}} f(w)  -u^\top w + \textstyle \frac{1}{2} \| w  \|_2^2.$ The dual of the problem 
can be derived in exactly the same way as shown in \eq{tvdual} in the previous section, using the definition of the support function  defined by 
\eq{suppTangentCone}. The corresponding dual is given by $ \max_{s \in \widehat{B}^{\mathcal{A}}(F)} -\textstyle \frac{1}{2} \| s - u\|_2^2, $ with the 
relationship $w=u-s$ at optimality. Thus, this corresponds to projecting $u$ onto the outer approximation of the base polytope, $\widehat{B}^{\mathcal{A}}(F)$,  which 
only has $m$ constraints instead of the $2^n-1$ constraints defining $B(F)$.  See an illustration in \myfig{outer}.

\subsection{Checking Optimality of a Basic Ordered Partition}
\label{sec:tvwarm}

Given a basic ordered partition $\mathcal{A}$, the associated $w \in \rb^n$ obtained from \eq{isotonic} is optimal for the TV problem in  \eq{tv} if and only if $s = u - w \in B(F)$ due to optimality conditions in \eq{tvdual}, which can be checked by minimizing the submodular function $F \!-\!s$. For a basic partition, a more efficient algorithm is available. 

By repeated application of submodularity, we have for all sets $C \subseteq V$, if $C_i = C \cap A_i$:

\BEAS
F(C) -s(C)
& = & F(V \cap C) - \sum_{i=1}^m s(C_i) \text{ (as $s$ is a modular function)}, \\
& = & F(B_m \cap C) -  \sum_{i=1}^m s(C_i)+ \sum_{i = 1}^{m-1} F( B_i \cap C) - F( B_i \cap C)  \text{ (as $B_m = V$)}, \\
& = & \sum_{i = 1}^m F(B_i \cap C) - F(B_{i-1} \cap C) - s(C_i)   \text{ (let $B_0 = \varnothing$ and as $F(\varnothing) = 0$)}, \\
& = & \sum_{i = 1}^m F((B_{i-1} \cup A_i) \cap C) - F(B_{i-1} \cap C) - s(C_i) 
\text{ (since $B_i = B_{i-1} \cup A_i$)},\\
& = & \sum_{i = 1}^m F((B_{i-1} \cap C) \cup (A_i \cap C)) - F(B_{i-1} \cap C) - s(C_i) ,\\
& = & \sum_{i = 1}^m F((B_{i-1} \cap C) \cup C_i) - F(B_{i-1} \cap C) - s(C_i) ,\\
& \geqslant &\sum_{i=1}^m \big[ F(B_{i-1} \cup C_i) - F(B_{i-1}) - s(C_i) \big] \\
&   & \hspace*{3cm} \text{  (as $(B_{i-1} \cap C) \subseteq B_{i-1}$ and due to submodularity of $F$). }
\EEAS
Moreover, we have  $s(A_i) = F(B_i) - F(B_{i-1})$, which implies $s(B_i) = F(B_i)$ for all $i \in \{1,\dots,m\}$, and thus all subproblems $
\min_{C_i \subseteq A_i} F(B_{i-1} \cup C_i) - F(B_{i-1}) - s(C_i)$
have non-positive values. This implies that we may check optimality by solving these $m$ subproblems: $s$ is optimal if and only if all of them have zero values.
This leads to smaller subproblems whose overall complexity is less than a single SFM oracle call. Moreover, for cut functions, it may be solved by a single oracle call on a graph where some edges have been removed~\citep{tarjan2006balancing}.

Given all sets $C_i$, we may then define a new ordered partition by splitting all $A_i$ for which $F(B_{i-1} \cup C_i) - F(B_{i-1}) - s(C_i)<0$. If no split is possible, the pair $(w,s)$ is optimal for \eq{tv}. Otherwise, this new strictly finer partition may not be basic, the value of the optimization problem in \eq{isotonic} is strictly lower  as shown in \mysec{convergence} (and leads to another basic ordered partition), which ensures the finite convergence of the algorithm. 

\subsection{Active-set Algorithm}
\label{sec:activesetAlg}
This leads to the novel active-set algorithm below.

\begin{algorithm}[H]
\KwData{Submodular function $F$ with SFM oracle, $u \in \rb^n$, ordered partition $\mathcal{A}$}
\KwResult{primal optimal: $w \in \rb^n$ and dual optimal: $s \in B(F)$}
\While {\texttt{True}}
{
Solve \eq{isotonic} using isotonic regression and update $\mathcal{A}$ with the basic ordered partition \;
Check optimality by solving $\min_{C_i \subseteq A_i} F(B_{i-1}  \cup C_i) \!-\! F(B_{i-1}) \!-\! s(C_i)$ for $i \in \{1,\dots,m\}$ \;
    \eIf{$s$ is optimal}
    {
        \texttt{break} \;
    }
    {
    for $i \in \{1, \dots, m\}$, split the set $A_i$ into $C_i$ and $A_i \setminus C_i$ in that order to get an updated ordered partition $\mathcal{A}$ \;
    }
}
\end{algorithm}

{\it Relationship with divide-and-conquer algorithm.}
 When starting from the trivial ordered partition $\mathcal{A} = (V)$, then  we exactly obtain a parallel version of the divide-and-conquer algorithm~\citep{groenevelt1991two}, that is, the isotonic regression problem is always solved without using the constraints of monotonicity, i.e., there are no merges, it is not necessary to re-solve the problems where nothing has changed. This shows that the number of iterations is then less than $n$. 

 The key added benefits in our formulation is the possibility of warm-starting, which can be very useful for building paths of solutions with different weights on the total variation. This is also useful for decomposable functions where many TV oracles are needed with close-by inputs. See experiments in \mysec{experiments}.

\subsection{Proof of Convergence}
\label{sec:convergence}
In order to prove the convergence of the algorithm in \mysec{activesetAlg}, we only need to show that if the optimality check fails in step (4), then step (7) introduces splits in the partition, which ensures that the isotonic regression in step (2) of the next iteration has a strictly lower value. Let us recall the isotonic regression problem solved in step (2):
\BEA
\min_{v \in \rb^m}
\sum_{i=1}^m \Big(v_i \big[ F(B_i) - F(B_{i-1}) - u(A_i) \big] + \textstyle \frac{1}{2} |A_i| v_i^2\Big)  \label{eq:isotonicCost} \\
                                \text{ such that }  v_1 \geqslant \cdots \geqslant v_m. \label{eq:isotonicConst}
\EEA
Steps (2) ensures that the ordered partition $\mathcal{A}$ is a basic ordered partition warranting that the inequality constraints are strict, 
i.e., no two partitions have the same value $v_i$ and the values $v_i$ for each element of the partition $i = \{1, \ldots, m\}$ is given through
\BEQ
\label{eq:isotonicOpt}
v_i |A_i| = u(A_i) - (F(B_i) - F(B_{i-1})),
\EEQ
which can be calculated in closed form.

The optimality check in step (4) decouples into checking the optimality in each subproblem as shown in \mysec{tvwarm}.
If the optimality test fails, then there is a subset of $C_i$ of $A_i$ for some of elements of the partition $\mathcal{A}$ such that 
\BEQ
F(B_{i-1} \cup C_i) - F(B_{i-1}) - s(C_i) < 0.
\label{eq:optimalFail}
\EEQ
We will show that the splits introduced by step (7) strictly reduces the function value of isotonic regression in \eq{isotonicCost}, while maintaining the 
feasibility of the problem. The splits modify the cost function of the isotonic regression as follows, as the objective function in \eq{isotonicCost} is equal to

\BEA
\sum_{i=1}^m \bigg( v_i\big[ F(B_{i-1} \cup C_i) -  F(B_{i-1}) - u(C_i)\big] + v_i \big[F(B_i) - F(B_{i-1} \cup C_i) - u(A_i \setminus C_i)\big]  \nonumber \\
 \ \ \ \ \ \ \ \ \ + \textstyle \frac{1}{2} v_i^2 |C_i|  + \textstyle \frac{1}{2} v_i^2 |A_i\setminus C_i|  \bigg) .\label{eq:Airestrict}
\EEA
Let us assume a positive $t \in \rb$, which is small enough. The direction that the isotonic regression moves is $v_i + t$ for the 
partition corresponding to $C_i$ and $v_i - t$ for the partition corresponding to $A_i \setminus C_i$ maintaining the feasibility of the 
isotonic regression problem, i.e.,  $v_1 \geqslant \cdots \geqslant v_i +t > v_i - t \geqslant \cdots \geqslant v_m$ . The function value is given by

\BEA
&   &\!\!\!\!\sum_{i=1}^m \Big( (v_i + t)\big[ F(B_{i-1} \cup C_i) -  F(B_{i-1}) - u(C_i)\big] \nonumber  \\
& & + (v_i - t)\big[F(B_i) - F(B_{i-1} \cup C_i) - u(A_i \setminus C_i)\big]  \nonumber + \textstyle \frac{1}{2}(v_i + t)^2 |C_i|  + \textstyle \frac{1}{2}(v_i - t)^2 |A_i \setminus C_i| \Big) \nonumber \\
& = &\sum_{i=1}^m \Big( \big( v_i\big[ F(B_{i-1} \cup C_i) -  F(B_{i-1}) - u(C_i)\big] + v_i \big[F(B_i) - F(B_{i-1} \cup C_i) - u(A_i \setminus C_i)\big]  \nonumber \\[-.1cm]
&   & \ + \textstyle \frac{1}{2} v_i^2 |C_i| + \textstyle \frac{1}{2}v_i^2|A_i \setminus C_i|  \big)  + \textstyle \frac{1}{2}t^2 |A_i|  \nonumber \\
&   & \ + t\  \big( 2F(B_{i-1} \cup C_i) - F(B_{i-1}) - F(B_{i})  - u(C_i) + u(A_i \setminus C_i) + v_i|C_i| - v_i|A_i\setminus C_i|\big)  \Big). \nonumber
\EEA
From this we can compute the directional derivative of the function at $t=0$, which is given by
\BEA
&   &2F(B_{i-1} \cup C_i) - F(B_{i-1}) - F(B_{i}) -u(C_i) + u(A_i \setminus C_i) + |C_i|v_i - |A_i\setminus C_i|v_i \nonumber \\
& = &2F(B_{i-1} \cup C_i) - F(B_{i-1}) - F(B_{i})  -2u(C_i) + u(A_i) + 2|C_i|v_i - |A_i|v_i \nonumber \\
& = &2\big(F(B_{i-1} \cup C_i) - F(B_{i-1}) - u(C_{i}) + v_i|C_i|\big)  \text{ (substituting \eq{isotonicOpt})} \nonumber\\
& = &2\big(F(B_{i-1} \cup C_i) - F(B_{i-1}) - s(C_{i})\big)  < 0 \text{ (as $s = u - w$ and \eq{optimalFail})}. \nonumber
\EEA
This shows that the function strictly decreases with the splits introduced in step (7).

\subsection{Discussion}
\label{sec:discuss}

{\it Certificates of optimality.}
The algorithm has dual-infeasible iterates $s$ (they only belong to $B(F)$ at convergence). Suppose that after step (3) we have $F(C) - s(C) \geqslant -\varepsilon$ for all $C \subset V$, where $\varepsilon$ shrinks as we run more iterations of the outer loop. This implies that
$s \in B(F + \varepsilon 1_{ {\rm Card} \in (1,n) })$, i.e., $s \in B(F_\varepsilon)$
with $F_\varepsilon = F + \varepsilon 1_{ {\rm Card} \in (1,n) }$. Since by construction $w = u - s$, we have:
\BEAS
f_\varepsilon(w) - u^\top w + \textstyle \frac{1}{2} \|w\|_2^2 + \textstyle \frac{1}{2} \| s  - u\|_2^2
& = &  \varepsilon \big| \max_{j \in V} w_j - \min_{j \in V} w_j \big|
 + f(w) - u^\top w + \| w\|^2 \\
& = &  \varepsilon \big| \max_{j \in V} w_j - \min_{j \in V} w_j \big| \\
& & \ \ + \sum_{i=1}^m v_i \big[ F(B_i) - F(B_{i-1}) - u(A_i) \big] +  \sum_{i=1}^m |A_i| v_i^2 \\
& = & \varepsilon \big| \max_{j \in V} w_j - \min_{j \in V} w_j \big| \text{ ( using \eq{isotonicOpt})}\\
& = & \varepsilon  \,  {\rm range}(w),
\EEAS
where  ${\rm range}(w) = \max_{k \in V} w_k - \min_{k \in V} w_k$. 
This means that $w$ is approximately optimal for $f(w) - u^\top w + \textstyle \frac{1}{2} \| w\|_2^2$ with \emph{certified gap} less than $   \varepsilon  \, {\rm range}(w) 
+ \varepsilon \,   {\rm range}(w^\ast)$.  

{\it Maximal range of an active-set solution.}
 For any ordered partition $\mathcal{A}$, and the optimal value of $w$ (which we know in closed form), we have
 $    {\rm range}(w) \leqslant     {\rm range}(u) + \max_{ i \in V} \big\{
 F(\{i\}) + F( V \backslash \{i\}) - F(V)  \big\}$. Indeed, for the $u$ part of the expression, this is because values of $w$ are averages of values of $u$; for the $F$ part of the expression, 
we always have by submodularity:
\BEAS
F(B_i) - F(B_{i-1}) & \leqslant & \sum_{k \in A_i} F(\{k\}) \text{ and} \\
F(B_i) - F(B_{i-1})&  \geqslant &- \sum_{k \in A_i} F(V) - F(V \backslash \{k\}).
\EEAS
This means that the certificate can be used in practice by replacing ${\rm range}(w^*)$ by its upperbound. See experimental evaluation for a 2D total variation denoising in \myapp{OptCert}.
 
{\it Exact solution.} If the submodular function only takes integer values and we have an approximate solution of the TV problem with gap $\varepsilon \leqslant \textstyle \frac{1}{4n}$, then we have the optimal solution~\citep{NIPS2014_5321}.
  
{\it Relationship with traditional active-set algorithm.} 
Given an ordered partition $\mathcal{A}$, an active-set method solves the unconstrained optimization problem in \eq{isotonic} to obtain a value of $v$ using the primary stationary conditions. The corresponding primal value $w = \sum_{i=1}^m v_i 1_{A_i}$ and dual value $s = u - w$ are optimal, if and only if,

\BEA
& \text{Primal feasibility :} \ \  w \in \mathcal{W}^{\mathcal{A}}, \label{eq:primalFeasibility} \\
& \text{Dual feasibility   :} \ \ s \in B(F). \label{eq:dualFeasibility}
\EEA

If \eq{primalFeasibility} is not satisfied, a move towards the optimal $w$ is performed to ensure primal feasibility by performing line search, i.e.,  two consecutive sets $A_i$ and $A_{i+1}$ with increasing values $v_i < v_{i+1}$ are merged and a potential $w$ is computed until primal feasibility is met. Then dual feasibility is checked and potential splits are proposed.

In our approach, we consider a different strategy which is more direct and does many merges simultaneously by using \emph{isotonic regression}. Our method explicitly moves from ordered partitions to ordered partitions and computes an optimal vector $w$, which is always feasible.

\section{Decomposable Problems}

\label{sec:decompalg}

Many interesting problems in signal processing and computer vision naturally involve submodular functions $F$ that decompose into $F = F_1 + \cdots + F_r$, with $r$ ``simple'' submodular functions~\citep{fot_submod}. For example, a cut function in a 2D grid decomposes into a function $F_1$ composed of cuts along vertical lines and a function $F_2$ composed of cuts along horizontal lines. For both of these functions, SFM oracles may be solved in $O(n)$ by message passing. For simplicity, in this paper, we consider the case $r=2$ functions, but following~\citet{komodakis2011mrf} and~\citet{treesubmod}, our framework easily extends to $r>2$.

\subsection{Reformulation as the Distance Between Two Polytopes}

Following~\citet{treesubmod}, we have the primal/dual problems :

\BEA
\min_{w \in \rb^n} f_1(w) + f_2(w)  -u^\top w + \textstyle \frac{1}{2} \| w  \|_2^2
\nonumber & = & \min_{w \in \rb^n}   \max_{s_1\in B(F_1), \ s_2 \in B(F_2) }  w^\top (s_1 + s_2) -u^\top w  + \textstyle \frac{1}{2} \| w  \|_2^2 \nonumber \\
& = &  \max_{s_1\in B(F_1), \ s_2 \in B(F_2) } \min_{w \in \rb^n}    (s_1+s_2 - u)^\top w + \textstyle \frac{1}{2} \| w  \|_2^2 \nonumber \\
& = &  \max_{s_1\in B(F_1), \ s_2 \in B(F_2) }  - \textstyle \frac{1}{2} \| s_1+s_2 -u  \|_2^2 ,
\label{eq:dualDecomp}
\EEA
with $w = u-s_1 - s_2$ at optimality.

This is the projection of $u$ on the sum of the base polytopes $B(F_1) + B(F_2) = B(F)$. Further, this may be interpreted as finding the distance between two polytopes $ B(F_1) - u/2$ and $ u/2 - B(F_2)$. Note that these two polytopes  typically do not intersect (they will if and only if $w=0$ is the optimal solution of the TV problem, which is an uninteresting situation). We now review Alternating projections (AP)~\citep{treesubmod}, Averaged alternating reflections (AAR)~\citep{treesubmod} and Dykstra's alternating projections (DAP)~\citep{accdyk} to show that a large number of total variation denoising problems need to be solved to obtain an optimal solution of \eq{dualDecomp}. The ability to warm start and solve these total variation denoising using our algorithm in \mysec{activesetAlg} can greatly improve the performance of each of these algorithms.

{\it Alternating projections (AP).}
The alternating projection algorithm~\citep{bauschke1997method} was proposed to solve the convex feasibility problem, i.e., to obtain a feasible point in the intersection of two polytopes. It is equivalent to performing block coordinate descent on the dual derived in \eq{dualDecomp}. Let us denote the projection onto a polytope $K$ as $\Pi_K$, i.e., $\Pi_K(y) = \argmax_{x \in K} - \| x - y \|_2^2.$
Therefore, alternating projections lead to the following updates for our problem.
$$
z_{t} = \Pi_{u/2 - B(F_2)}\Pi_{B(F_1) - u/2}(z_{t-1}),
$$
where $z_{0}$ is an arbitrary starting point. Thus each of these steps require TV oracles for $F_1$ and $F_2$ since projection onto the base polytope is equivalent to performing TV denoising as shown in \eq{tvdual}.

{\it Averaged alternating reflections (AAR).}
The averaged alternating reflection algorithm~\citep{bauschke2004finding}, which is also known as Douglas-Rachford splitting can be used to solve convex feasibility problems. It is observed to converge quicker than alternating projection~\citep{treesubmod,seshTV} in practice. We now introduce a reflection operator for the polytope $K$ as $R_K$, i.e., $R_K = 2\Pi_K - I,$ where $I$ is the identity operator. Therefore, reflection of $t$ on a polytope $K$ is given by  $R_K(t) = 2\Pi_K(t) - t$. The updates of each iteration of the averaged alternating reflections, which starts with an auxiliary sequence $z_0$ initialized to $0$ vector, are given by 
$$
z_{t} = \textstyle \frac{1}{2}(I + R_{u/2 - B(F_2)}R_{B(F_1) - u/2})(z_{t-1}).
$$
In the feasible case, i.e., intersecting polytopes, the sequence~$z_t$ weakly converges to a point in the intersection of the polytopes. However, in our case, we have non intersecting polytopes which leads to a converging sequence of $z_t$ with AP but a diverging sequence of $z_t$ with AAR. However, when we project $z_t$ by using the projection operation, $s_{1,t}  =  \Pi_{B(F_1) - u/2}(z_{t}); s_{2,t}  =  \Pi_{u/2 - B(F_2)}(s_{1,t})$ the sequences $s_{1,t}$ and $s_{2,t}$ converge to the nearest points on the polytopes, $B(F_1) - u/2$ and $u/2 - B(F_2)$~\citep{bauschke2004finding} respectively.

{\it Dykstra's alternating projections (DAP).} Dykstra's alternating projection algorithm~\citep{bauschke1994dykstra} retrieves a convex feasible point closest to an arbitrary point, which we assume to be $0$. It can also be used and has a form of primal descent interpretation, i.e., as coordinate descent for a well-formulated primal problem~\citep{gaffke1989cyclic}. Let us denote $\iota_K$ as the indicator function of a convex set $K$.  In our case we consider finding the nearest points on the polytopes $B(F_1) - u/2$ and $u/2 - B(F_2)$ closest to $0$, which can be formally written as:

\BEAS
&   & \min_{\substack{s \in B(F_1) - \frac{u}{2} \\  s \in \frac{u}{2} - B(F_2)}} \frac{1}{2} \| s \|^2_2 \\
& = & \min_{s \in \rb^n} \frac{1}{2} \| s \|^2_2  + \iota_{B(F_1) - \frac{u}{2}}(s) + \iota_{\frac{u}{2} - B(F_2)}(s) \\
& = & \min_{s \in \rb^n} \frac{1}{2} \| s \|^2_2  + \iota_{B(F_1) - \frac{u}{2}}(s) + \iota_{B(F_2) - \frac{u}{2}}(-s) \\
& = & \min_{s \in \rb^n} \bigg( \frac{1}{2} \| s \|^2_2  + \max_{w_1 \in \rb^n}  w_1^{\top}s - f_1(w_1) + \frac{w_1^{\top}u}{2} + \max_{w_2 \in \rb^n}  - w_2^{\top}s - f_2(w_2) + \frac{w_2^{\top}u}{2}\bigg)  \\
& = & \max_{\substack{w_1 \in \rb^n\\ w_2 \in \rb^n}}  \bigg(- f_1(w_1) - f_2(w_2) + \frac{(w_1 + w_2)^{\top}u}{2} + \min_{s \in \rb^n} \bigg(\frac{1}{2} \| s \|^2_2  + (w_1 - w_2)^{\top}s \bigg) \bigg) \\
& = & \max_{\substack{w_1 \in \rb^n\\ w_2 \in \rb^n}}  - f_1(w_1) - f_2(w_2) + \frac{(w_1 + w_2)^{\top}u}{2}  - \frac{1}{2} \| w_1 - w_2\|^2_2 \\
& = & \min_{\substack{w_1 \in \rb^n\\ w_2 \in \rb^n}}  f_1(w_1) + f_2(w_2) - \frac{(w_1 + w_2)^{\top}u}{2}  + \frac{1}{2} \| w_1 - w_2\|^2_2,
\EEAS 
where $s = w_2 - w_1$ at optimality. The block coordinate descent algorithm then gives

\BEAS
w_{1,t} & = & \text{prox}_{f_1 - u/2}(w_{2, t-1}) =  (I - \Pi_{B(F_1) - u/2})(w_{2, t-1}), \\
s_{1,t} & = & w_{2, t-1} - w_{1, t}, \\
w_{2,t} & = & \text{prox}_{f_2 - u/2}(w_{1, t}) = (I - \Pi_{B(F_2) - u/2})(w_{1, t}), \\
s_{2,t} & = & w_{1, t} - w_{2, t},
\EEAS
where $I$ is the identity matrix. This is exactly the same as Dykstra's alternating projection steps. 

\begin{figure}[t]
\begin{center}
\begin{tabular}{cc}
\includegraphics[width=0.45\textwidth]{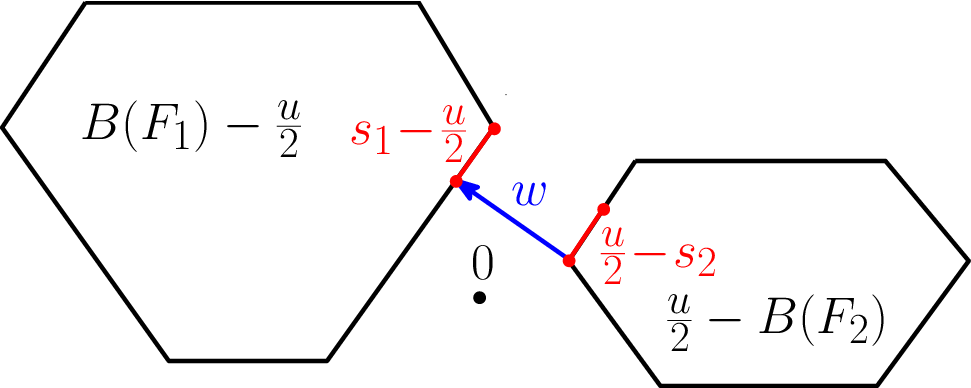} & 
\includegraphics[width=0.45\textwidth]{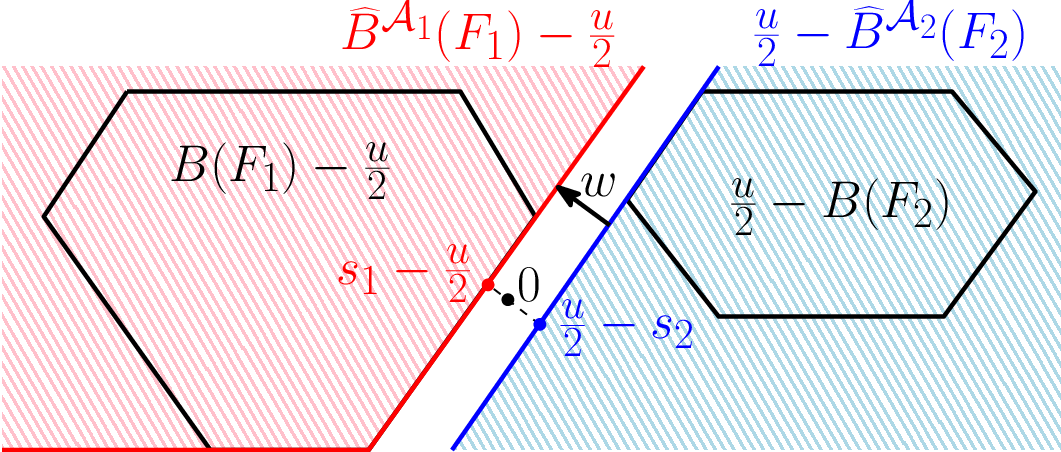}  \\
(a) & (b) 
\end{tabular}
\end{center}

\vspace*{-.5cm}

\caption{ {Closest point between two polytopes.} (a) Output of Dykstra's alternating projection algorithm for the TV problem, the pair $(s_1,s_2)$ may not be unique while $w = s_1 +s_2 - u$ is. (b) Dykstra's alternating projection output for outer approximations.}
\label{fig:outer-closest}
\end{figure}

We have implemented it, and it behaves similar to alternating projections, but it still requires TV oracles for projection (see experiments in \mysec{experiments}). There is however a key difference: while alternating projections and alternating reflections always converge to a pair of closest points, Dykstra's alternating projection algorithm converges to a \emph{specific} pair of points, namely the pair closest to the initialization of the algorithm~\citep{bauschke1994dykstra}; see an illustration in \myfig{outer-closest}-(a). This insight will be key in our algorithm to avoid cycling.  

Assuming TV oracles are available for $F_1$ and $F_2$, \citet{treesubmod} and \citet{seshTV} use alternating projection~\citep{bauschke1997method} and alternating reflection~\citep{bauschke2004finding} algorithms to solve dual optimization problem in \eq{dualDecomp} . \citet{nishihara2014convergence} gave a extensive theoretical analysis of alternating projection and showed that it converges linearly. However, these algorithms are equivalent to block \emph{dual} coordinate descent and cannot be cast explicitly as descent algorithms for the primal TV problem. On the other hand, Dykstra's alternating projection is a descent algorithm on the primal, which enables local search over partitions. Complex TV oracles are often implemented by using SFM oracles recursively with the divide-and-conquer strategy on the individual functions. Using our algorithm in \mysec{activesetAlg}, they can be made more efficient using warm-starts (see experiments in \mysec{experiments}).

\subsection{Local Search over Partitions using Active-set Method}
Given our algorithm for a single function, it is natural to perform a local search over two partitions $\mathcal{A}_1$ and $\mathcal{A}_2$, one for each function $F_1$ and $F_2$, and consider in the primal formulation a  weight vector~$w$ compatible with both $\mathcal{A}_1$ and $\mathcal{A}_2$; or, equivalently, in the dual formulation, two outer approximations
$\widehat{B}^{\mathcal{A}_1}(F_1)$ and $\widehat{B}^{\mathcal{A}_2}(F_2)$. That is, given the ordered partitions  $\mathcal{A}_1$ and $\mathcal{A}_2$, using a similar derivation as in \eq{dualDecomp}, we obtain the primal/dual pairs of optimization problems
\BEQ
\max_{\substack{ s_1 \in \widehat{B}^{\mathcal{A}_1}(F_1) \\  s_2 \in \widehat{B}^{\mathcal{A}_2}(F_2)}} \textstyle  -\frac{1}{2} \| u - s_1 - s_2 \|_2^2   = \min_{\substack{ w \in \mathcal{W}^{\mathcal{A}_1} \\ w \in \mathcal{W}^{\mathcal{A}_2} }} f_1(w) + f_2(w) - u^\top w +\textstyle  \frac{1}{2} \| w\|_2^2, \nonumber
\EEQ
with $w = u - s_1 - s_2$ at optimality. 

{\it Primal solution by isotonic regression.}
The primal solution~$w$ is unique by strong convexity. Moreover, it has to be compatible with both 
$\mathcal{A}_1$ and $\mathcal{A}_2$, which is equivalent to being compatible with the \emph{coalesced} ordered partition $\mathcal{A} = {\rm coalesce}(\mathcal{A}_1,\mathcal{A}_2)$ defined as the coarsest ordered partition compatible by both. As shown in \myapp{coalescing}, $\mathcal{A}$ may be found in time $O(\min(m_1, m_2)n)$.  

Given $\mathcal{A}$, the primal solution $w$ of the subproblem may be found by isotonic regression like in \mysec{isotonic} in time $O(m)$ where $m$ is the number of sets in $\mathcal{A}$. However, finding the optimal dual variables $s_1$ and $s_2$ turns out to be more problematic. We know that $s_1 + s_2 = u - w$ and that $ s_1 + s_2 \in \widehat{B}^{\mathcal{A}}(F)$, but the split of $s_1 + s_2$ into $(s_1, s_2)$ is unknown.

{Obtaining dual solutions.} Given ordered partitions $\mathcal{A}_1$ and $\mathcal{A}_2$, a unique well-defined pair $(s_1, s_2)$ could be obtained by using  convex feasibility algorithms
such as alternating projections \citep{bauschke1997method} or alternating reflections~\citep{bauschke2004finding}. However, the result would depend in non understood ways on  the 
initialization, and we have observed cycling of the active-set algorithm. Using Dykstra's alternating projection algorithm allows us to converge to a unique well-defined pair 
$(s_1,s_2)$ that will lead to a provably non-cycling algorithm.  

When running the Dykstra's alternating projection algorithm starting from $0$ on the polytopes $\widehat{B}^{\mathcal{A}_1}(F_1) - u/2$ and $u/2 - \widehat{B}^{\mathcal{A}_2}(F_2)$, if $w$ is the unique distance vector between the two polytopes, then the iterates converge to the projection of $0$ onto the convex sets of elements in the two polytopes that achieve the minimum distance~\citep{bauschke1994dykstra}. See \myfig{outer-closest}-(b) for an illustration. This algorithm is however slow to converge when the polytopes do not intersect. Note that $w \neq 0$ in most of our situations and convergence is hard to monitor because primal iterates of the Dykstra's alternating projection diverge~\citep{bauschke1994dykstra}.
 
\begin{figure}[t]
\begin{center}
\begin{tabular}{cc}
\includegraphics[width=0.45\textwidth]{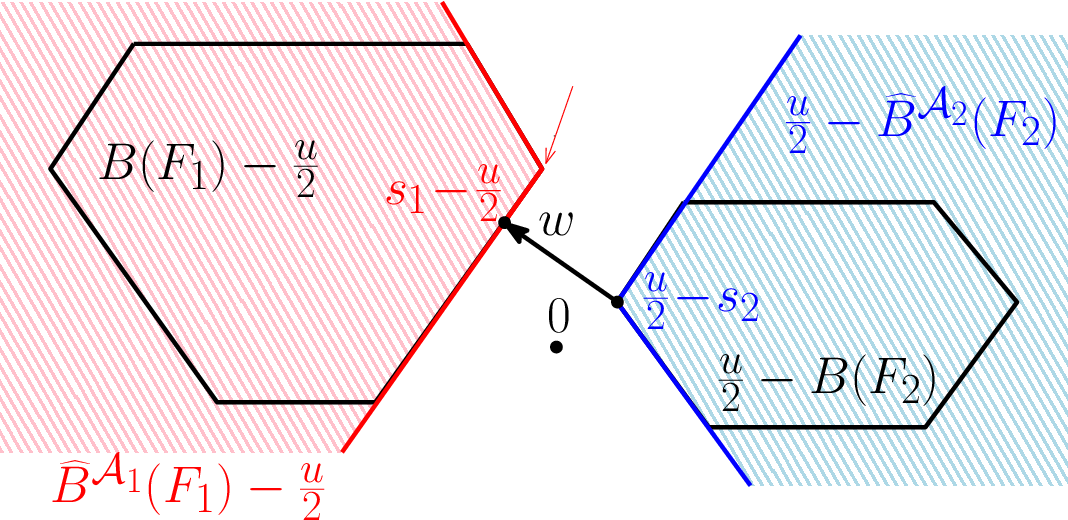}  & 
\includegraphics[width=0.45\textwidth]{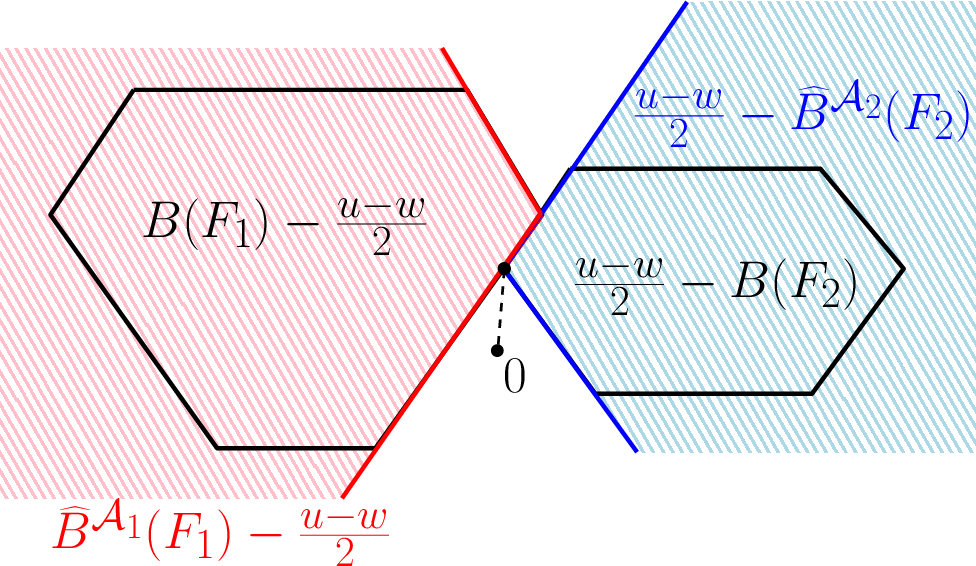}  \\
(a) & (b)
\end{tabular}
\end{center}

\vspace*{-.5cm}

\caption{ {Translated intersecting polytopes.} (a) output of our algorithm before translation. (b)~Translated formulation.}
\label{fig:outer-trans}
\end{figure}

{Translated intersecting polytopes.} In our situation, we have more to work with than just the ordered partitions: \emph{we also know the vector $w$} (as mentioned earlier, it is obtained  cheaply from isotonic regression). Indeed, from Lemma 2.2 and Theorem 3.8  from~\citet{bauschke1994dykstra}, given this vector $w$, we may translate the two polytopes and now obtain a formulation where the two polytopes do intersect; that is we aim at projecting~$0$ on the (non-empty) intersection of $\widehat{B}^{\mathcal{A}_1}(F_1) - u/2 + w/2$ and $u/2-w/2-\widehat{B}^{\mathcal{A}_2}(F_2)$. See \myfig{outer-trans}. We also refer to this as the {\em translated Dykstra problem}\footnote{We refer to finding a Dykstra solution for translated intersecting polytopes as translated Dykstra problem.} in the rest of the paper. This is equivalent to solving the following optimization problem

\BEA
& &  \min_{\substack{s \in \widehat{B}^{\mathcal{A}_1}(F_1) - \frac{u - w}{2} \\  s \in \frac{u - w}{2} - \widehat{B}^{\mathcal{A}_2}(F_2)}} \textstyle \frac{1}{2} \| s \|^2_2 \\
& = & \min_{s \in \rb^n} \textstyle \frac{1}{2} \| s \|^2_2  + \iota_{\widehat{B}^{\mathcal{A}_1}(F_1) - \frac{u - w}{2}}(s) + \iota_{\frac{u - w}{2} - \widehat{B}^{\mathcal{A}_2}(F_2)}(s),\nonumber  \\
& = & \min_{s \in \rb^n} \textstyle \frac{1}{2} \| s \|^2_2  + \iota_{\widehat{B}^{\mathcal{A}_1}(F_1) - \frac{u - w}{2}}(s) + \iota_{\widehat{B}^{\mathcal{A}_2}(F_2) - \frac{u - w}{2} }(-s),\nonumber  \\
& = & \min_{s \in \rb^n} \bigg(\textstyle \frac{1}{2} \| s \|^2_2  + \max_{w_1 \in \mathcal{W}^{\mathcal{A}_1}} w_1^{\top}s - f_1(w_1) + \frac{w_1^{\top}(u - w)}{2}  \nonumber \\[-0.8em]
&   & \hspace*{3cm} + \max_{w_2 \in \mathcal{W}^{\mathcal{A}_2}} - w_2^{\top}s - f_2(w_2) + \frac{w_2^{\top}(u - w)}{2} \bigg), \nonumber \\
& = & \max_{\substack{w_1 \in \mathcal{W}^{\mathcal{A}_1}\\ w_2 \in \mathcal{W}^{\mathcal{A}_2}}}  \bigg(- f_1(w_1) - f_2(w_2) + \frac{(w_1 + w_2)^{\top}(u - w)}{2} + \min_{s \in \rb^n} \frac{1}{2} \| s \|^2_2  + (w_1 - w_2)^{\top}s \bigg), \nonumber \\
& = & \max_{\substack{w_1 \in \mathcal{W}^{\mathcal{A}_1}\\ w_2 \in \mathcal{W}^{\mathcal{A}_2}}}\bigg(  - f_1(w_1) - f_2(w_2) + \frac{(w_1 + w_2)^{\top}(u - w)}{2}
- \frac{1}{2} \| w_1 - w_2\|^2_2\bigg), \nonumber \\
& = & \!\!\!\! \min_{\substack{w_1 \in \mathcal{W}^{\mathcal{A}_1}\\ w_2 \in \mathcal{W}^{\mathcal{A}_2}}} \bigg( f_1(w_1) + f_2(w_2) - \frac{(w_1 + w_2)^{\top}(u - w)}{2}  + \frac{1}{2} \| w_1 - w_2\|^2_2\bigg), \label{eq:dykPrimal}
\EEA 
with $s = w_2 - w_1$ at optimality.

In \mysec{linear} we propose algorithms  to solve the above optimization problems. Assuming that we are able to solve this step efficiently, we now present our active-set algorithm for decomposable functions below.

\subsection{Active-set Algorithm for Decomposable Functions}
\label{sec:alg}

\begin{algorithm}[H]
\KwData{Submodular functions $F_1$ and $F_2$ with SFM oracles, $u \in \rb^n$, ordered partitions $\mathcal{A}_1, \mathcal{A}_2$}
\KwResult{primal optimal: $w \in \rb^n$ and dual optimal: $s_1 \in B(F_1), s_2 \in B(F_2)$}
\While {\texttt{True}}
{
$\mathcal{A} = {\rm coalesce}(\mathcal{A}_1,\mathcal{A}_2)$ \;
Estimate $w$ by solving $\min_{w \in \mathcal{W}^{\mathcal{A}}} f(w) - u^\top w + \frac{1}{2} \|w\|_2^2$ using isotonic regression \;
{\bf Projection step:} Estimate $s_1 \in \widehat{B}^{\mathcal{A}_1}(F_1)$ and $s_2 \in \widehat{B}^{\mathcal{A}_2}(F_2)$ by projecting $0$ onto the intersection of $\widehat{B}^{\mathcal{A}_1}(F_1) - u/2 + w/2$ and $u/2-w/2-\widehat{B}^{\mathcal{A}_2}(F_2)$ using any of the algorithms described in \mysec{linear} \;
Merge the sets in $\mathcal{A}_j$ which are tight for $s_j$, $j \in \{1,2\}$\;
Check optimality of $s_1$ and $s_2$ as described in \mysec{tvwarm}\;
    \eIf{$s_1$ and $s_2$ are optimal}
    {
        \texttt{break} \;
    }
    {
    for $j \in \{1, 2\}$ and $i_j \in \{1, \dots, m_j\}$, split the set $A_{j, i_j}$ into $C_{j, i_j}$ and $A_{j, i_j} \setminus C_{j, i_j}$ in that order to get an updated ordered partition $\mathcal{A}_j$ \;
    }
}
\end{algorithm}
 
Given two ordered partitions $\mathcal{A}_1$ and $\mathcal{A}_2$, we obtain $s_1 \in \hat{B}^{\mathcal{A}_1}(F_1)$ and $s_2 \in \hat{B}^{\mathcal{A}_2}(F_2)$ as described in the following section. 
The solution $w = u - s_1 - s_2$ is optimal if and only if both $s_1 \in B(F_1)$ and $s_2 \in B(F_2)$. When checking the  optimality described in \mysec{tvwarm}, we split the partition. As shown in \myapp{optimalitydec}, either (a) $\|w\|_2^2$ strictly increases at each iteration, or (b) $\|w\|_2^2$ remains constant but $\|s_1-s_2\|_2^2$ strictly increases. This implies that the algorithm is finitely convergent.

\subsection{Optimizing the ``Translated Dykstra Problem"}
\label{sec:linear}
In this section, we describe algorithms for the ``projection step" of the active-set algorithm proposed in \mysec{alg} that optimizes the translated Dykstra problem in \eq{dykPrimal}, i.e., 
\BEQ
\min_{\substack{w_1 \in \mathcal{W}^{\mathcal{A}_1} \\  w_2 \in \mathcal{W}^{\mathcal{A}_2}}} f_1(w_1) + f_2(w_2) - \textstyle \frac{(w_1 + w_2)^{\top}(u-w)}{2} + \textstyle \frac{1}{2} \| w_1 - w_2 \|^2.
\label{eq:primalDyk1}
\EEQ
The corresponding dual optimization problem is given by 
\BEQ
\min_{\substack{s \in \widehat{B}^{\mathcal{A}_1}(F_1) - \frac{u - w}{2} \\  s \in \frac{u - w}{2} - \widehat{B}^{\mathcal{A}_2}(F_2)}} \textstyle \frac{1}{2} \| s \|^2_2,
\label{eq:dualDyk}
\EEQ
with the optimality condition $s = w_2 - w_1$. Note that the only link to submodularity is that $f_1$ and $f_2$ are linear functions on $\mathcal{W}^{\mathcal{A}_1}$ and $\mathcal{W}^{\mathcal{A}_2}$, respectively. The rest of this section primarily deals with optimizing a quadratic program and we present two algorithms in  \mysec{Dykstra} and \mysec{pas}.

\subsubsection{Accelerated Dykstra's Algorithm}
\label{sec:Dykstra}

In this section, we find the projection of the origin onto the intersection of the translated base polytopes obtained by solving the optimization problem in \eq{dykPrimal} given by 
\BEQ
\min_{\substack{w_1 \in \mathcal{W}^{\mathcal{A}_1} \\  w_2 \in \mathcal{W}^{\mathcal{A}_2}}} f_1(w_1) + f_2(w_2) - \textstyle \frac{(w_1 + w_2)^{\top}(u-w)}{2} + \textstyle \frac{1}{2} \| w_1 - w_2 \|^2,
\nonumber
\EEQ
using Dykstra's alternating projection. It can be solved using the following Dykstra's iterations:
 \BEAS
s_{1,t} & = &  \Pi_{\widehat{B}^{\mathcal{A}_1} (F_1)} (u/2 - w/2 + w_{2,t-1}   ), \\ 
w_{1,t} & = &  u/2 - w/2 + w_{2,t-1}  - s_{1,t}  ,  \\
s_{2,t} & = &  \Pi_{\widehat{B}^{\mathcal{A}_2}(F_2)}  ( u/2 - w/2 + w_{1,t} ), \\  
w_{2,t} & = &  u/2 - w/2 + w_{1,t} - s_{2,t},
\EEAS
with $\Pi_C$ denoting the orthogonal projection onto the sets $C$, solved here by isotonic regression.  Note that the value of the auxiliary variable $w_{2}$ can be warm-started. The algorithm converges linearly for polyhedral sets~\citep{shusheng2000estimation}.

In our simulations, we have used the recent accelerated version of~\citet{accdyk}, which led to faster convergence. In order to monitor convergence of the algorithm, we compute the value of $\|u-w-s_{1,t}-s_{2,t}\|_1$, which is equal to zero at convergence. Note that the algorithm is not finitely convergent and gives only {\em approximate} solutions. Therefore, we introduce the approximation parameters $\varepsilon_1$ and $\varepsilon_2$ such that $s_1$ lies in the $\epsilon_1$-neighborhood of the base polytope of $F_1$, i.e., $\min_{A \subseteq V} F_1(A) - s_1(A) \geq -\varepsilon_1$ and $s_2$ lies in the $\varepsilon_2$-neighborhood of the base polytope of $F_2$, i.e., $\min_{A \subseteq V} F_2(A) - s_2(A) \geq -\varepsilon_2$, respectively. See \myapp{alpha} for more details on the approximation. The optimization problem can also be decoupled into smaller optimization problems by using the knowledge of the face of the base polytopes on which $s_1$ and $s_2$ lie. This is still slow to converge in practice and therefore we present an active-set method in the next section.

\subsubsection{Primal Active-set Method}
\label{sec:pas}
In this section, we find the projection of the origin onto the intersection of the translated base polytopes given by \eq{dykPrimal} using the standard active-set method~\citep{Nocedal2006} by solving a set of linear equations. For this purpose, we derive the equivalent optimization problems using equality constraints. 

The ordered partition, $\mathcal{A}_j$ is given by $(A_{j, 1}, \ldots, A_{j, m_j})$, where $m_j$ is the number of elements in the ordered partitions. Let $B_{j, i_j}$ be defined as $(A_{j, 1} \cup \cdots \cup A_{j, i_j})$. Therefore,

\BEA
f_j(w_j) & = &  \sum_{i_j = 1}^{m_j} v_{j, i_j} \bigg(F_j(B_{j, i_j}) -  F_j(B_{j, i_j-1})\bigg) \label{eq:fval} \\
w_j      & = &  \sum_{i_j = 1}^{m_j} v_{j, i_j} 1_{A_{j, i_j}} \label{eq:wval} \\
& & \text{with the constraints, } v_{j,1} \geqslant  \cdots \geqslant v_{j,m_j}. \label{eq:constraints}
\EEA

On substituting \eq{fval}, \eq{wval} and \eq{constraints} in \eq{dykPrimal}, we have an equivalent optimization problem:
\BEAS
\min_{\substack{v_{1,1} \geq  \ldots \geq v_{1,m_1}\\  v_{2,1} \geq  \ldots \geq v_{2,m_2}}} & \displaystyle\sum_{i_1 = 1}^{m_1} & \bigg(F_1(B_{1, i_1})  -  F_1(B_{1, i_1-1})  - \frac{u(A_{1, i_1}) - w(A_{1, i_1})}{2}\bigg)v_{1, i_1}  \\ 
+ & \displaystyle\sum_{i_2 = 1}^{m_2} & \bigg(F_2(B_{2, i_2}) -  F_2(B_{2, i_2-1})  - \frac{u(A_{2, i_2}) - w(A_{2, i_2})}{2}\bigg)v_{2, i_2} \\ 
+ & \displaystyle\sum_{i_1 = 1}^{m_1} & \frac{1}{2} |A_{1, i_1}| v_{1, i_1}^2+ \sum_{i_2 = 1}^{m_2} \frac{1}{2} |A_{2, i_2}| v_{2, i_2}^2  \displaystyle\sum_{i_1 = 1}^{m_1} \sum_{i_2 = 1}^{m_2} v_{1, i_1} v_{2, i_2} 1_{A_{1, i_1}}^{\top} 1_{A_{2, i_2}}.
\EEAS

This can be written as a quadratic program in $ x = \left( \begin{array}{c} v_1 \\ v_2 \end{array} \right) $ with inequality constraints in the following form

\BEQ
\min_{\substack{x \in \rb^{m_1 + m_2}\\ D({\mathcal{A}_1, \mathcal{A}_2})x \succcurlyeq 0}} \frac{1}{2} x^{\top}Q({\mathcal{A}_1, \mathcal{A}_2})x + c({\mathcal{A}_1, \mathcal{A}_2})^{\top}x.
\label{eq:qp}
\EEQ
Here, $D(\mathcal{A}_1, \mathcal{A}_2)$ is a sparse matrix of size $(m_1 + m_2 -2) \times (m_1+m_2)$, which is a block diagonal matrix
containing the difference or first order derivative matrices of sizes $m_1 - 1 \times m_1$ and $m_2 - 1 \times m_2$ as the blocks 
and $c({\mathcal{A}_1, \mathcal{A}_2})$ is a linear vector that can be computed using the function evaluations of $F_1$ and $F_2$. 
Note that these evaluations need to be done only once.

\begin{figure}
\begin{center}
\includegraphics[width=0.45\textwidth]{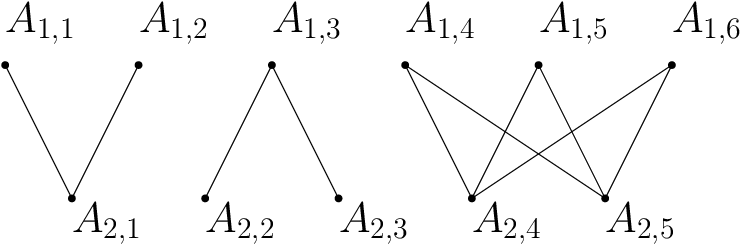}  
\end{center}

\vspace*{-.5cm}

\caption{Bipartite graph to compute $Q({\mathcal{A}_1, \mathcal{A}_2})$ with $\mathcal{A}_1$ having $m_1 = 6$ components and $\mathcal{A}_2$ having $m_2 = 5$.}
\label{fig:bipartite}
\end{figure}

{Computing $Q({\mathcal{A}_1, \mathcal{A}_2})$.} Let us consider a bipartite graph, $G = (\mathcal{A}_1, \mathcal{A}_2, E)$, with $m_1 + m_2$ nodes 
representing the ordered partitions of $\mathcal{A}_1$ and $\mathcal{A}_2$ respectively. The weight of the edge between each 
element of ordered partitions of $\mathcal{A}_1$, represented by $A_{1, i_1}$ and each element of ordered partitions of $\mathcal{A}_2$, 
represented by $A_{2, i_2}$ is the number of elements of the ground set $V$ that lie in both these partitions and can be written as
$ e(A_{1, i_1}, A_{2, i_2}) = 1_{A_{1, i_1}}^{\top} 1_{A_{2, i_2}}$ for all $e \in E$. The matrix $Q({\mathcal{A}_1, \mathcal{A}_2})$ 
represents the Laplacian matrix of the graph $G$. \myfig{bipartite} shows a sample bipartite graph with $m_1 = 6$ and $m_2 = 5$.

{Initializing with primal feasible point.} Primal active-set methods start with a primal feasible point and continue to maintain
primal feasible iterates. In our case, the starting point may be obtained using the weight vector $w$ that is estimated using isotonic regression. 
The vector $w$ is compatible both with $\mathcal{A}_1$ and $\mathcal{A}_2$, i.e., $w \in \mathcal{W}^{\mathcal{A}_1}$ and $w \in \mathcal{W}^{\mathcal{A}_2}$. 
Therefore, we may obtain the vectors $v_1$ and $v_2$ from $w$ and initialize the primal feasible starting point using $v_1$ and $v_2$.

Optimizing the quadratic program in \eq{qp} by using active-set methods is equivalent to finding the face of the constraint set on which 
the optimal solution lies. For this purpose, we need to be able to solve the quadratic program in \eq{qp} with equality constraints.

{Equality constrained QP.} Let us now consider the following quadratic program with equality constraints
\BEQ
\min_{\substack{p \in \rb^{m_1 + m_2}\\ D'p = 0}} \frac{1}{2} p^{\top}Q({\mathcal{A}_1, \mathcal{A}_2})p + \big(Q({\mathcal{A}_1, \mathcal{A}_2})x_k + c({\mathcal{A}_1, \mathcal{A}_2})\big)^{\top}p,
\label{eq:eqqp}
\EEQ
where $D'$ is the subset of the constraints in $D({\mathcal{A}_1, \mathcal{A}_2})$, i.e., indices of constraints that are tight and $x_k$ is a primal-feasible point. We refer the indices of the tight constraints as the {\em working set} and represent them by the set $W$ in the algorithm. Therefore, the set of constraints in $D'$ is the restriction of the constraint set $D({\mathcal{A}_1, \mathcal{A}_2})$ to the working set constraints denoted by $W$. The vector $p$ gives the direction of strict descent of the cost function in \eq{qp} from feasible point $x_k$~\citep{Nocedal2006}.

Without loss of generality, let us assume that the equality constraints are $v_{j, k_j} = v_{j, k_j + 1}$ for any $k_j$ in $[0, m_j)$. Let $\mathcal{A}'_j$ be the new ordered partition formed 
by merging $A_{j, k_j}$ and $A_{j, k_j + 1}$ as $v_{j, k_j} = v_{j, k_j + 1}$. Similarly, $x'_t$ can be computed from $x_t$ by merging the weights $v_{j, k_j}$ and $v_{j, k_j + 1}$ into a single weight for the merged element of the ordered partition. Finding the optimal vector $p'$ using the quadratic program in \eq{eqqp} with respect to the ordered partition $\mathcal{A}'_j$ is equivalent to solving the following unconstrained quadratic problem,

\BEA
\mathcal{Q}(\mathcal{A}'_1, \mathcal{A}'_2, x'_t)  =  \min_{p' \in \rb^{m'_1 + m'_2}} \bigg(\frac{1}{2} p'^{\top}Q({\mathcal{A}'_1, \mathcal{A}'_2})p'+ \big(Q({\mathcal{A}'_1, \mathcal{A}'_2})x'_t + c({\mathcal{A}'_1, \mathcal{A}'_2})\big)^{\top}p'\bigg),
\label{eq:unqp}
\EEA
where $m'_j$ is the number of elements of the ordered partition ${\mathcal{A}'}_j$. This can be estimated by solving a linear system using conjugate gradient descent. The complexity of each iteration of the conjugate gradient is given by $O((m'_1 + m'_2)k)$ where $k$ is the number of non-zero elements in the sparse matrix, $Q({\mathcal{A}'_1, \mathcal{A}'_2})$~\citep{vishnoi2013lx}. We can build $p$ from $p'$ by repeating the values for the elements of the partition that were merged.

{Primal active-set algorithm.} We now can describe the standard primal active-set method.

\begin{algorithm}[H]
\KwData{Laplacian matrix, $Q({\mathcal{A}_1, \mathcal{A}_2})$ and vector, $c({\mathcal{A}_1, \mathcal{A}_2})$, ordered partitions $\mathcal{A}_1, \mathcal{A}_2$}
\KwResult{$x^* \in \rb^{m_1 + m_2}$}
\textbf{Initialize}: 
$t=0$\;
Primal feasible $x_0$ using the solution of isotonic regression $w$\;
$W_0$ with indices of rows of $D(\mathcal{A}_1, \mathcal{A}_2)$ that are equal to 0\;
\While {\texttt{True}}
{
    \textbf{Estimate primal}: Solve equality constrained QP in \eq{eqqp} with equality constraints indexed by working set $W_t$ to find optimal $p$ \;
    \eIf{ p == 0 }  
    {
        \textbf{Estimate dual}: $\lambda = D(\mathcal{A}_1, \mathcal{A}_2)^{- \top} \big(Q({\mathcal{A}_1, \mathcal{A}_2})x_t + c({\mathcal{A}_1, \mathcal{A}_2})\big)$\; 
        \eIf{ $\lambda_i \geq 0$ for all $i \in W_t$}
        { 
            \texttt{break} \;
        }
        {
            update $j = \argmin_{j \in W_t} \lambda_j$ \;
            update $W_{t+1} = W_t \setminus \{j\}$\;
            update $x_{t+1} = x_t$ \;
        }
    }
    {
        \textbf{Line search}: Find least $\beta$ that retains feasibility of $x_{t+1} = x_t + \beta p$ and find the blocking constraints $B_t$ \;
        $W_{t+1} = W_t \cup B_t$ \;
    }
    $t = t+1$
}
\Return {$x^* = x_t$}
\end{algorithm}
 
We can estimate $w_1$ and $w_2$ from $x^*$, which will enable us to estimate $s$ feasible in \eq{dualDyk}. Therefore we can estimate the dual 
variable $s_1 \in B^{\mathcal{A}_1}(F_1)$ and $s_2 \in B^{\mathcal{A}_2}(F_2)$ using~$s$.

\subsection{Decoupled Problem}
 In our context, the quadratic program in \eq{qp} can be decoupled into smaller optimization problems. Let us consider the bipartite graph $G = (\mathcal{A}_1, \mathcal{A}_2, E)$ of which $Q$ is the Laplacian matrix. The number of connected components of the graph, $G$, is equal to the number of level-sets of $w$. 

Let $m$ be the total number of connected components in~$G$. These connected components define a partition on the 
ground set $V$ and a total order on elements of the partition can be obtained using the levels sets of $w$. Let $k$ 
denote the index of each bipartite subgraph of $G$ represented by $G_k = (\mathcal{A}_{1,k}, \mathcal{A}_{2, k}, E_k)$,
where $k = 1, 2, \ldots, m$. Let $J_k$ denote the indices of the nodes of $G_k$ in $G$.
\BEQ
x^*_{J_k} = \argmin_{\substack{x \in \rb^{m_{1,k} + m_{2,k}}\\ D(\mathcal{A}_1, \mathcal{A}_2)_kx \succcurlyeq 0}} \frac{1}{2} x^{\top}Q(\mathcal{A}_1, \mathcal{A}_2)_{J_kJ_k}x + c(\mathcal{A}_1, \mathcal{A}_2)_{J_k}^{\top}x, 
\label{eq:decoupqp}
\EEQ
where $m_{j,k}$ is size of $\mathcal{A}_{j,k}$. Therefore, $m_{1,k} + m_{2,k}$ is the total number of nodes in the subgraph $G_k$.
Note that this is exactly equivalent to decomposition of the base polytope of $F_j$ into base polytopes of submodular functions formed
by contracting $F_j$ on each individual component representing the connected component $k$. See  \myapp{decouple} for more details.

\section{Experiments}
\label{sec:experiments}

In this section, we show the results of the algorithms proposed on various problems. We first consider the problem of solving total variation denoising for a non decomposable function using
active-set methods in \mysec{expNonDecomp}, specifically cut functions. Here, our experiments mainly focus on the time comparisons with state-of-art methods and also show an important setting where
we show the gain due to the ability to warm-start our algorithm. In \mysec{expDecomp}, we consider cut functions on a 3D grid decomposed into a function of the 2D grid and a function of chains. We then consider a 2D grid and a concave function on cardinality, which is not a cut function. Our algorithm leads to marginal gains for the usual non decomposable functions. However, in the non decomposable case there are many total variation problems to be solved. The ability to warm-start lends to huge improvements when compared to the usage of standard total variation oracles in this setting.
\begin{figure*}
\begin{center}
\begin{tabular}{cc}
\includegraphics[width=0.45\textwidth]{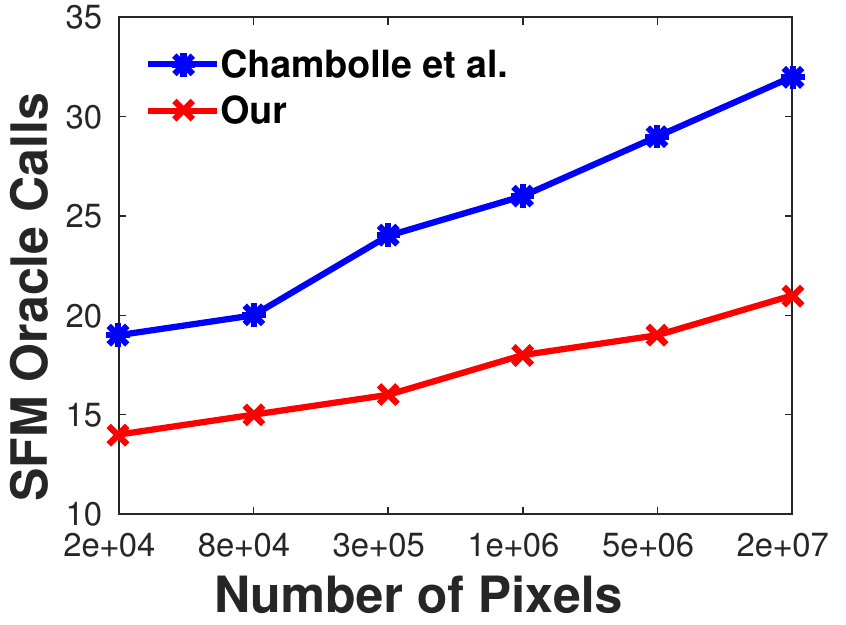} &
\includegraphics[width=0.45\textwidth]{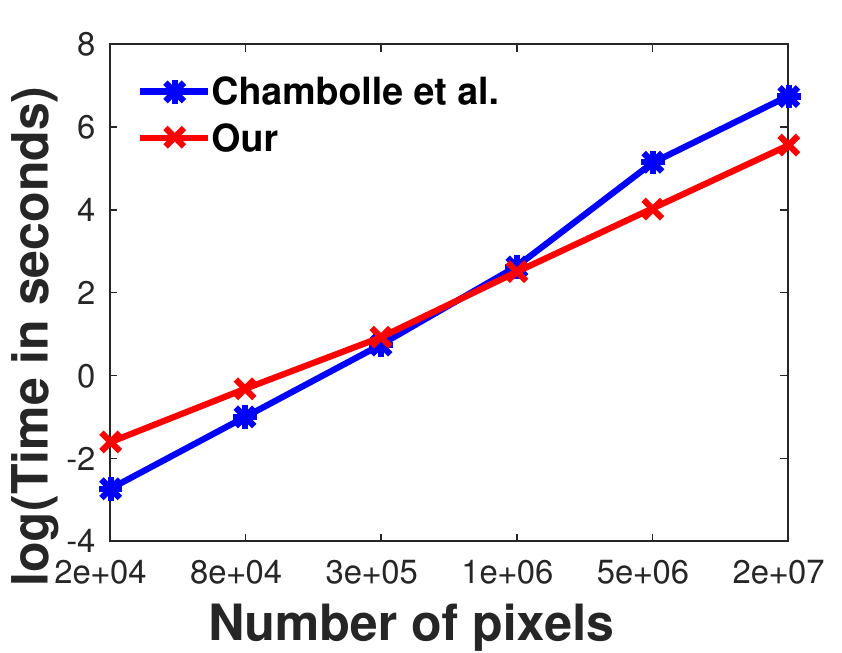} \\
(a) & (b) \\
\includegraphics[width=0.45\textwidth]{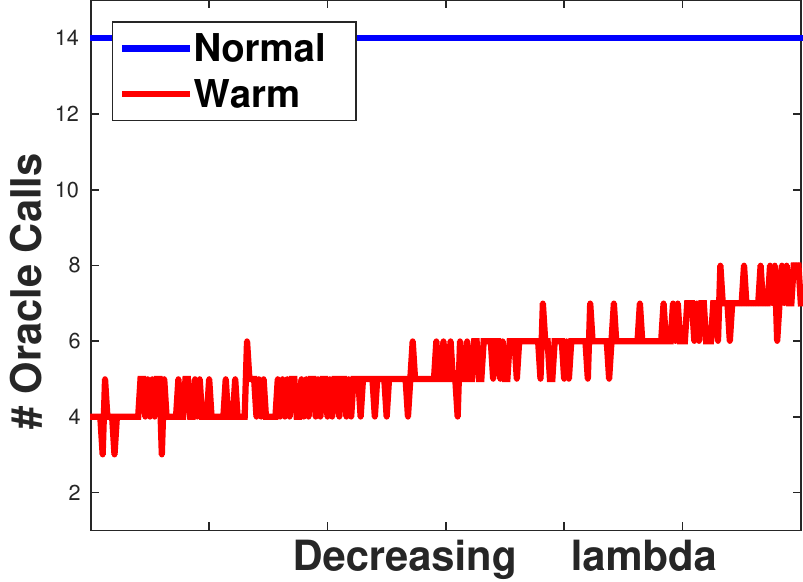} &
\includegraphics[width=0.45\textwidth]{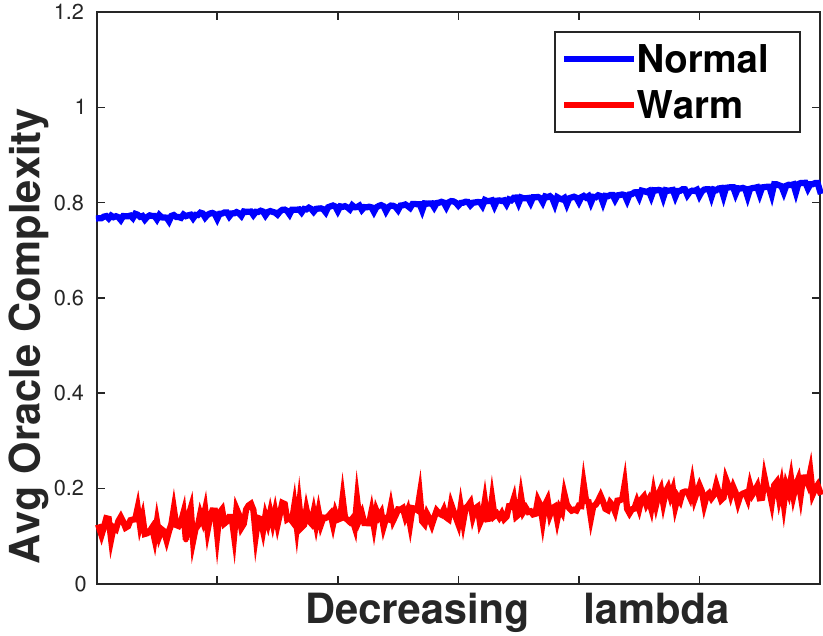} \\
(c) & (d)
\end{tabular}
\end{center}

\vspace*{-.5cm}

\caption{(a) Number of SFM oracle calls for images of various sizes, (b) Time taken for images of various sizes, (c) Number of iterations with and without warm start, (d) Average complexity of the oracle with and without warm start. }
\label{fig:2DSFM}
\end{figure*}

\subsection{Non-decomposable Total Variation Denoising}
\label{sec:expNonDecomp}

Our experiments consider images, which are 2-dimensional grids with vertex neighborhood of size~4. The data set comprises of 6 different images of varying sizes. We consider a large image of size $5616 \times 3744$ and recursively scale into a smaller image of half the width and half the height maintaining the aspect ratio. Therefore, the size of each image is four times smaller than the size of the previous image. We restrict to anisotropic uniform-weighted total variation to compare with~\citet{chambolle2009total} but our algorithms work as well with weighted total variation, which is
standard in computer vision, and on any graph with SFM oracles. Therefore, the unweighted total variation is 
$$f(w)~=~\lambda~\sum_{i \sim j}~|w_i~-~w_j|,$$ 
where $\lambda$ is a regularizing constant for solving the total variation problem in \eq{tv}. 

Maxflow~\citep{boykov04} is used as the SFM oracle for checking the optimality of the ordered partitions. \myfig{2DSFM}(a) shows the number of SFM oracle calls required to solve the TV problem for images of various sizes. Note that for the algorithm of \citet{chambolle2009total} each SFM oracle call optimizes smaller problems sequentially, while each SFM oracle call in our method optimizes several independent smaller problems in parallel. Therefore, our method has lesser number of oracle calls than~\citet{chambolle2009total}. However, oracle complexity of each call is higher for the our method when compared to~\citet{chambolle2009total}. \myfig{2DSFM}(b) shows the time required for each of the methods to solve the TV problem to convergence. We have an optimized code and only use the oracle as plugin which takes about 80-85 percent of the running time. This is primarily the reason our approach takes more time than \citep{chambolle2009total} despite having fewer oracle calls for small images. 

\myfig{2DSFM}(c) also shows the ability to warm start by using the output of a related problem, i.e., when computing the solution for several values of~$\lambda$ (which is typical in practice). In this case, we use optimal ordered partitions of the problem with larger~$\lambda$ to warm start the problem with smaller~$\lambda$. It can be observed that warm start of the algorithm requires lesser number of oracle calls to converge than using  the initialization with trivial ordered partition. Warm start also largely helps in reducing the burden on the SFM oracle. With warm starts the number of ordered partitions does not change much over iterations.  Hence, it suffices to query only ordered partitions that have changed. To analyze this we define {\em oracle complexity} as the fraction of pixels in the elements of the partitions that need to be queried.  Oracle complexity is averaged over iterations to understand the average burden on the oracle per iteration. With warm starts this reduces drastically, which can be observed in \myfig{2DSFM}(d). 
\subsection{Decomposable Total Variation Denoising and SFM}
\label{sec:expDecomp}
\begin{figure*}
\begin{center}
\begin{tabular}{ccc}
\includegraphics[width=0.3\textwidth]{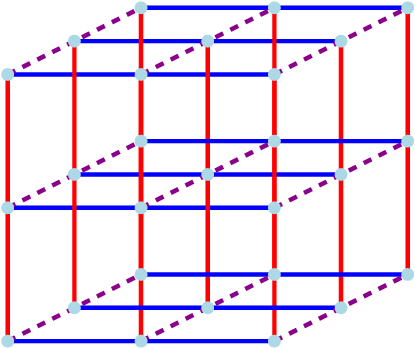}&
\includegraphics[width=0.3\textwidth]{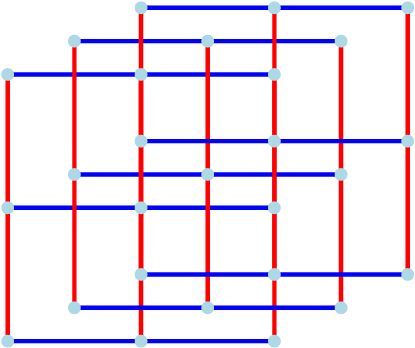}&
\includegraphics[width=0.3\textwidth]{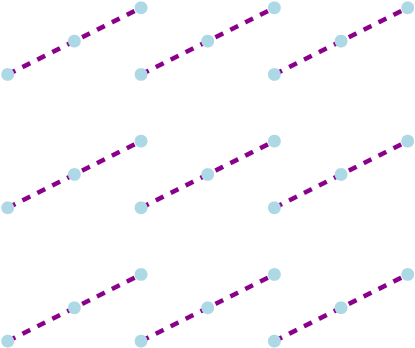}\\
(a) & (b) & (c)
\end{tabular}
\end{center}

\vspace*{-.5cm}

\caption{(a) Cut function $F$ defined on a 3D-grid may be decomposed into: (b) $F_1$ represented by solid edges (red and blue) and (c) $F_2$ represented by dashed lines (magenta)}
\label{fig:3Dgrid}
\end{figure*}

{\it Cut functions. }  In the decomposable case, we consider the SFM and TV problems on a cut function defined on a 3D-grid. The 3D-grid consists of lines parallel lines in each dimension as shown in \myfig{3Dgrid}. It can be decomposed into two functions $F_1$ and $F_2$, where $F_1$ is composed on parallel 2D-grids and $F_2$ is composed of parallel chains. From \myfig{3Dgrid}, the function $F_1$ represents all the solid edges(red and blue) whereas the function $F_2$ represents the dashed edges(magenta). For brevity, we refer to each 2D-grid of the function $F_1$ as a {\em frame}. The SFM oracle for the function $F_1$ is the maxflow-mincut~\citep{boykov04} algorithm, which may run in parallel for all frames. Similarly the SFM oracle for the function $F_2$ is the message passing algorithm, which may run in parallel for all chains. The corresponding TV oracles, i.e., {\em projection algorithm} for $F_1$ and $F_2$ may be solved using the algorithm described in \mysec{activesetAlg} due to availability of the respective SFM oracles. We consider averaged alternating reflection (AAR)~\citep{bauschke2004finding} by solving each projection without warm-start and counting the total number of SFM oracle calls of $F_1$ to solve the SFM and TV on the 3D-grid as our baseline. (SGD-P) denotes the dual subgradient based method~\citep{komodakis2011mrf} modified with Polyak's rule~\citep{poljak1987introduction} to solve SFM on the 3D-grid. We show the performance of alternating projection (AP-WS), averaged alternating reflection (AAR-WS)~\citep{bauschke2004finding} and Dykstra's alternating projection (DAP-WS)~\citep{bauschke1994dykstra} using warm start of each projection with the ordered partitions. WS denotes warm start variant of each of the algorithm. The performance of the active-set algorithm proposed in \mysec{alg} with inner loop solved using the primal active-set method proposed in \mysec{pas} is represented by (ACTIVE).

\begin{figure*}
\begin{center}
\begin{tabular}{cc}
\includegraphics[width=0.45\textwidth]{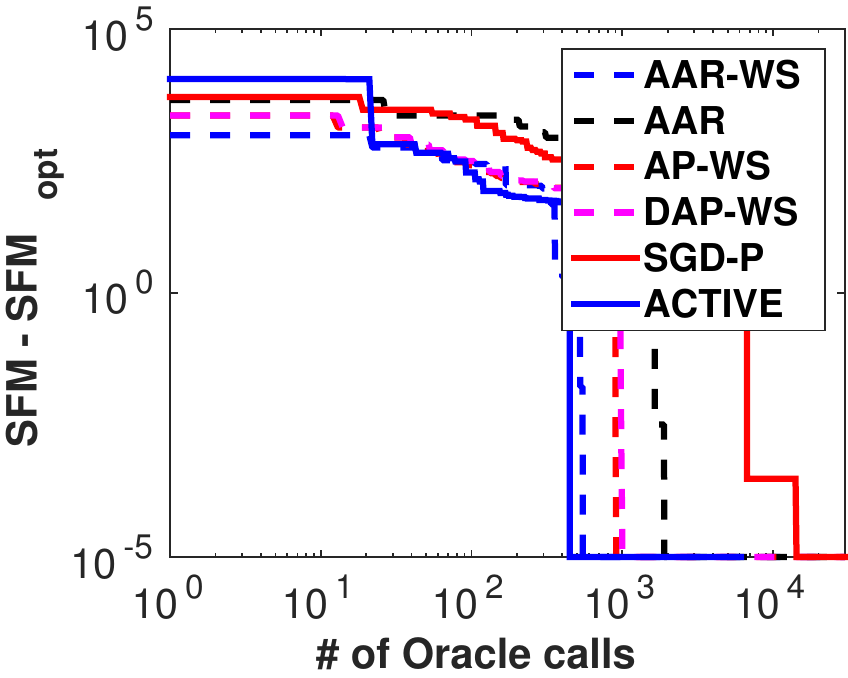} &
\includegraphics[width=0.45\textwidth]{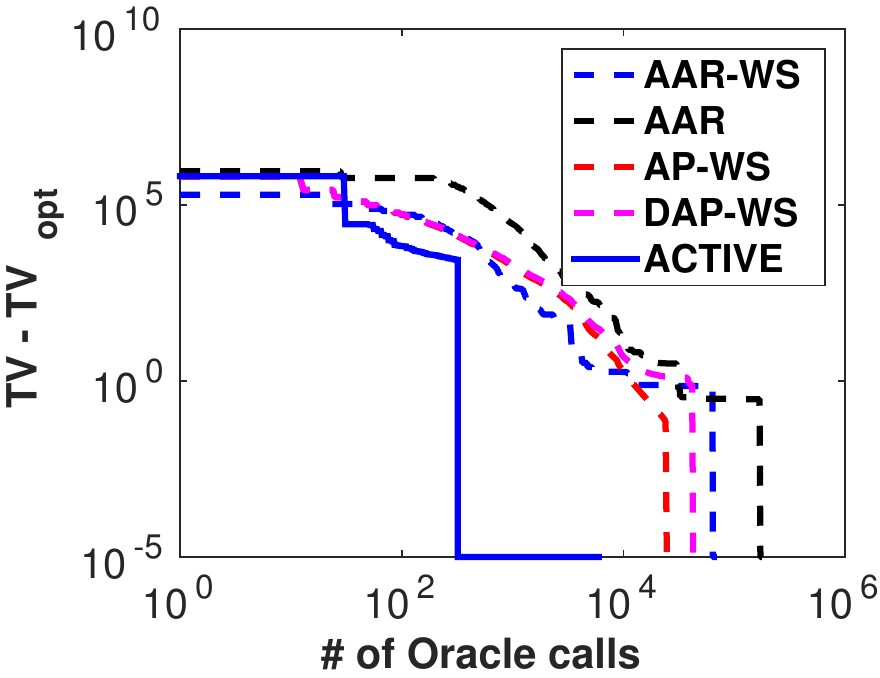} \\
(a) & (b) \\
\includegraphics[width=0.45\textwidth]{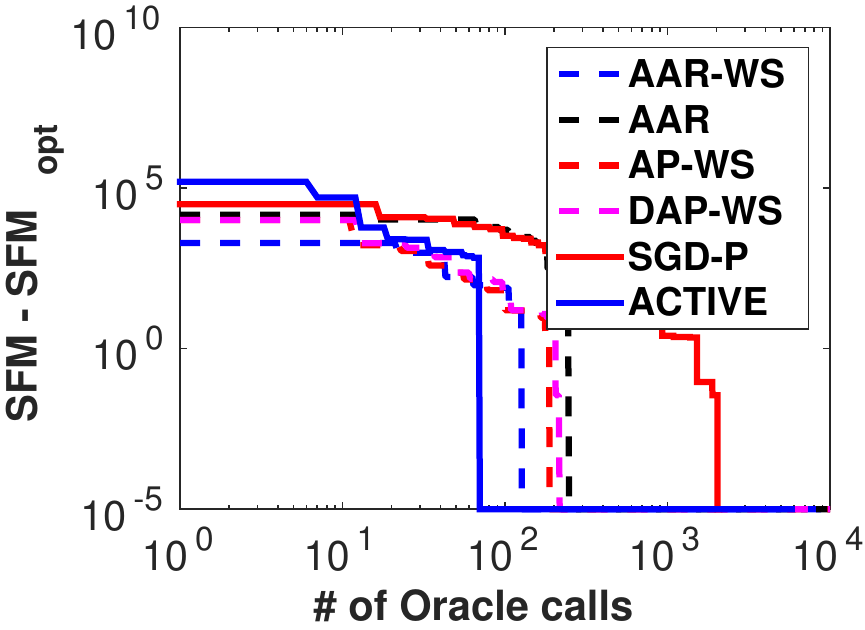} &
\includegraphics[width=0.45\textwidth]{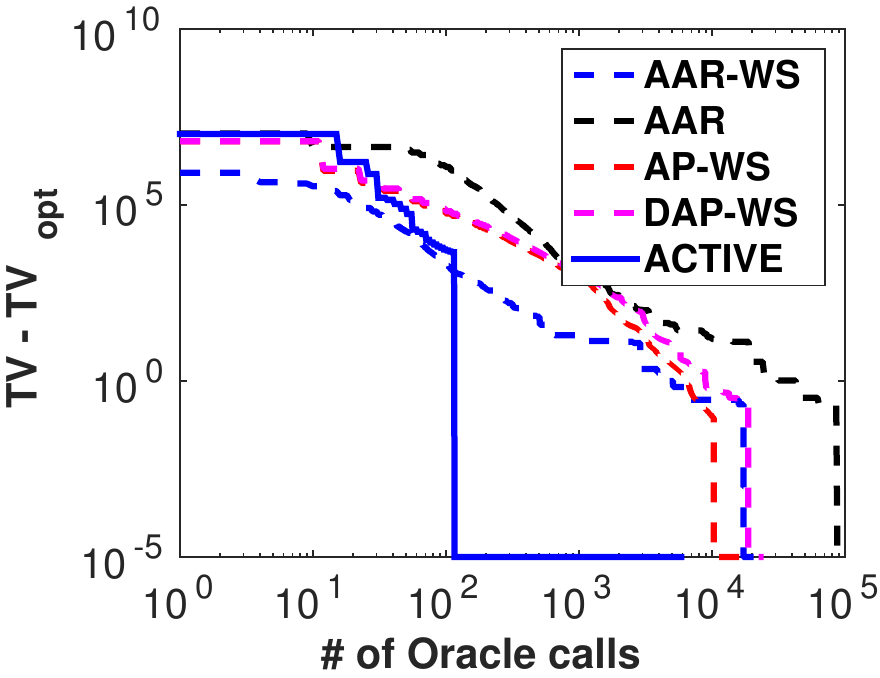} \\
(c) & (d)
\end{tabular}
\end{center}

\vspace*{-.5cm}

\caption{(a) Number of 2D SFM calls to obtain 3D SFM, (b) Number of 2D SFM calls to obtain 3D TV, (c) Number of 2D SFM calls to obtain SFM of 2D + concave function, (d) Number of 2D SFM calls to obtain TV of 2D + concave function. }
\label{fig:expDecomp}
\end{figure*}

In our experiments, we consider the 3D volumetric data set of the Stanford bunny~\citep{maxflowDataset} of size $102 \times 100 \times 79$.  The function $F_1$ represents $102$ frames while $F_2$ represents the $7900$ chains. The dimension of each frame in $F_1$ is $100 \times 79$, while the length of each chain in $F_2$ is $102$. \myfig{expDecomp}~(a) and (b) show that (AP-WS), (AAR-WS), (DAP-WS) and (ACTIVE) require relatively less number of oracle calls when compared to when compared to AAR or SGD-P. Note that we count 2D SFM oracle calls as they are more expensive than the SFM oracles on chains. 

{\it Time comparisons. } We also performed time comparisons between the iterative methods and the combinatorial methods on standard data sets. The standard mincut-maxflow~\citep{boykov04} on the 3D volumetric data set of the Standard bunny~\citep{maxflowDataset} of size $102 \times 100 \times 79$ takes $0.11$ seconds while averaged alternating reflections (AAR) without warm start takes $0.38$ seconds. The averaged alternating reflections with warm start (AAR-WS) takes $0.21$ seconds and the active-set method (ACTIVE) takes $0.38$ seconds. The main bottleneck in the active-set method is the inversion of the Laplacian matrix and it could considerably improve by using methods suggested by~\citet{vishnoi2013lx}. Note that the projection on the base polytopes of $F_1$ and $F_2$ can be parallelized by projecting onto each of the 2D frame of $F_1$ and each line of $F_2$ respectively~\citep{seshTV}. The times for (AAR), (AAR-WS) and (ACTIVE) use parallel multi-core architectures\footnote{20 core, Intel(R) Xeon(R) CPU E5-2670 v2 @ 2.50GHz with 100Gigabytes of memory. We only use up to 16 cores of the machine to ensure accurate timings.} while the combinatorial algorithm only uses a single core. Note that cut functions on grid structures are only a small subclass of submodular functions with such efficient combinatorial algorithms. In contrast, our algorithm works on more general class of sum of submodular functions than just with cut functions.

{\it Concave functions on cardinality. }  In this experiment we consider our SFM problem of sum of a 2D cut on a graph of size $5616 \times 3744$ and a super pixel based concave function on cardinality~\citep{stobbe11,treesubmod}. The unary potentials of each pixel is calculated using the Gaussian mixture model of the color features. The edge weight $a(i,j) = \exp(- \| y_i - y_j\|^2)$, where $y_i$ denotes the RGB values of the pixel $i$. In order to evaluate the concave function, regions $R_j$ are extracted via superpixels and, for each $R_j$, defining the function $F_2(S) = |S||R_j \setminus S|$. We use 200 and 500 regions. \myfig{expDecomp}~(c) and (d) shows that (AP-WS), (AAR-WS), (DAP-WS) and (ACTIVE) algorithms converge for solving TV quickly by using only SFM oracles and relatively less number of oracle calls. Note that we count 2D SFM oracle calls.

\section{Conclusion}
In this paper, we present an efficient active-set algorithm for optimizing quadratic losses regularized by \lova extension of a submodular function using the SFM oracle of the function. We also present an active-set algorithms to minimize sum of ``simple'' submodular functions using SFM oracles of the individual ``simple'' functions. We also show that these algorithms are competitive to the existing state-of-art algorithms to minimize submodular functions.

\acks{We acknowledge support from the European Research Council grant SIERRA (project 239993). K.~S.~Sesh~Kumar also acknowledges the support from the European Research Council grant of Prof.~Vladimir~Kolmogorov DOiCV(project 616160) at IST Austria. The comments of the reviewers have helped us improve the presentation significantly.}

\appendix

\section{Certificates of Optimality}
\label{app:OptCert}
We consider a total variation denoising problem on an 2D image of dimensions $384 \times 288$ using the algorithm proposed in \mysec{activesetAlg} for non-decomposable functions, where we assume a 2D SFM oracle. Here, we plot the optimality gap given by $(f(w) - u^{\top}w + \textstyle \frac{1}{2}\|w\|^2_2) - (f(w^\ast) - u^{\top}w^\ast + \textstyle \frac{1}{2}\|w^\ast\|^2_2)$ and the bound, $\varepsilon  \, {\rm range}(w) + \varepsilon \,   {\rm range}(w^\ast)$, proposed in \mysec{discuss}. Here, $w$ is the solution at the end of each iteration of the algorithm and $w^\ast$ is the optimal solution.
\begin{figure}
\begin{center}
\includegraphics[width=0.5\textwidth]{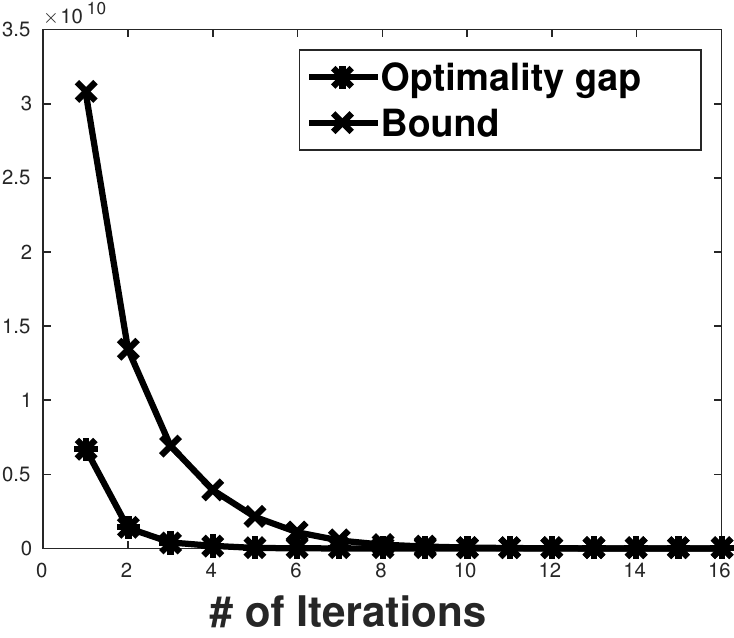}
\end{center}

\vspace*{-.5cm}

\caption{Certified gap vs. bound for TV denoising on a 2D image using the algorithm in \mysec{activesetAlg}.}
\label{fig:gap}
\end{figure}

\section{Algorithms for Coalescing Partitions}
\label{app:coalescing}

The basic interpretation in coalescing two ordered partitions is as follows. Given an ordered partition $\mathcal{A}_1$ and $\mathcal{A}_2$ with $m_1$ and $m_2$ elements in the partitions respectively, 
we define for each $j=1,2$,  $\forall i_j = (1,\ldots,m_j)$, $B_{j,i_j} = (A_{j,1} \cup \ldots \cup A_{j,i_j})$. 
The inequalities defining the outer approximation of the base polytopes are given by
hyperplanes defined by $$\forall i_j = (1, \ldots, m_j), s_j(B_{j, i_j}) \leq F_j(B_{j, i_j}).$$ The hyperplanes defined by common sets of both these partitions, defines the coalesced ordered partitions.
The following algorithm performs coalescing between these partitions.
\BIT
\item \textbf{Input}: Ordered partitions $\mathcal{A}_1$ and $\mathcal{A}_2$. 
\item \textbf{Initialize}: $x=1$, $y=1$, $z= 1$ and $C = \varnothing$.
\item \textbf{Algorithm}: Iterate until $x=m_1$ and $y=m_2$
\BNUM
 \item If $|B_{1,x}| > |B_{2,y}|$ then $y := y + 1$.
 \item If $|B_{1,x}| < |B_{2,y}|$ then $x := x + 1$.
 \item If $B_{1,x} == B_{2,y}$ then
 \BIT
     \item $A_{z} = (B_{1,x} \setminus C)$, 
     \item $C= B_{1,x}$, and 
     \item $z:=z+1$.
 \EIT
\ENUM
\item \textbf{Output}: $m = z$, ordered partitions $\mathcal{A} = (A_1, \ldots, A_m)$.
\EIT

{Running time.} The algorithm terminates in $\min(m_1, m_2)$ iterations and the checking condition for step (3) takes $n$ iterations. Therefore, the algorithm overall takes a time of $O(\min(m_1, m_2) n)$.

\section{Optimality of Algorithm for Decomposable Problems}
\label{app:optimalitydec}

In step (10) of the algorithms, when we split partitions, the value of  the primal/dual pair of optimization algorithms
\BEQ
\max_{\substack{ s_1 \in \widehat{B}^{\mathcal{A}_1}(F_1) \\ s_2 \in \widehat{B}^{\mathcal{A}_2}(F_2)}} \textstyle  -\frac{1}{2} \| u - s_1 - s_2 \|_2^2 =   \min_{\substack{ w \in \mathcal{W}^{\mathcal{A}_1}\\ w \in \mathcal{W}^{\mathcal{A}_2}} } f_1(w) + f_2(w) - u^\top w +\textstyle  \frac{1}{2} \| w\|_2^2, \nonumber
\EEQ
cannot increase. This is because, when splitting, the constraint set for the minimization problem only gets bigger. Since at optimality, we have $w = u - s_1 - s_2$, $\| w\|_2$ cannot decrease, which shows the first statement.

Now, if $\|w\|_2$ remains constant after an iteration, then it has to be the same (and not only have the same norm), because  the optimal $s_1$ and $s_2$ can only move in the direction orthogonal to $w$.

In step (4) of the algorithm,   we  project $0$ on the (non-empty) intersection of $\widehat{B}^{\mathcal{A}_1}(F_1) - u/2 + w/2$ and $u/2-w/2-\widehat{B}^{\mathcal{A}_2}(F_2)$. This corresponds to minimizing $\frac{1}{2} \| s_1 - u/2 + w/2 \|^2$ such that $s_1 \in \widehat{B}^{\mathcal{A}_1}(F_1)$ and $s_2 = u - w - s_1 \in \widehat{B}^{\mathcal{A}_2}(F_2)$. This is equivalent to minimizing $\frac{1}{8} \| s_1 - s_2\|^2$. We have: 
 \BEAS
 \max_{\substack{ s_1 \in \widehat{B}^{\mathcal{A}_1}(F_1)  \\ s_2 \in \widehat{B}^{\mathcal{A}_2}(F_2)\\ s_1 + s_2 = u - w}} - \frac{1}{8} \| s_1 - s_2 \|_2^2 & = & 
  \min_{\substack{w_1 \in \mathcal{W}^{\mathcal{A}_1} \\ w_2 \in \mathcal{W}^{\mathcal{A}_2}}}   \max_{\substack{s_1 \in \rb^n \\  s_2 \in \rb^n \\ s_1 + s_2 = u - w}}
   \bigg(- \frac{1}{8} \| s_1 - s_2 \|_2^2 + f_1(w_1) + f_2(w_2) \\ [-2em]
  &   & \ \ \ \ \ \ \ \ \ \ \ \ \ \ \ \ \ \ \ 
        \ \ \ \ \ \ \ \ \ \ \ \ \ \ \ \ \ \ \ - w_1^\top s_1 - w_2^\top s_2 \bigg)\\
  & = & 
  \min_{\substack{w_1 \in \mathcal{W}^{\mathcal{A}_1} \\ w_2 \in \mathcal{W}^{\mathcal{A}_2}}}   \max_{s_2 \in \rb^n}
  \bigg( - \frac{1}{8} \|u - w -  2 s_2 \|_2^2 + f_1(w_1) + f_2(w_2) \\ [-2em]
  &   & \ \ \ \ \ \ \ \ \ \ \ \ \ \ \ \ \ \ \ 
        \ \ \ \ \ \ \ \ \ \ \ \ \ \ \ \ \ \ \ - w_1^\top ( u - w - s_2) - w_2^\top s_2\bigg)   \\
  & = & 
  \min_{\substack{w_1 \in \mathcal{W}^{\mathcal{A}_1} \\ w_2 \in \mathcal{W}^{\mathcal{A}_2}}}   
\max_{  
  s_2 \in \rb^n}
  \bigg(- \frac{1}{8} \|u - w \|_2^2  -  \frac{1}{2} \| s_2 \|_2^2 
  + \frac{1}{2} s_2^\top ( u - w ) \\ [-1.5em]
 &   &  \ \ \ \ \ \ \ \ 
        \ \ \ \ \ \ \ \ + f_1(w_1) + f_2(w_2) - w_1^\top ( u - w - s_2) - w_2^\top s_2  \bigg) \\
  & = & 
  \min_{\substack{w_1 \in \mathcal{W}^{\mathcal{A}_1} \\ w_2 \in \mathcal{W}^{\mathcal{A}_2}}}   
  - w_1^\top ( u - w )  + f_1(w_1) + f_2(w_2)
  - \frac{1}{8} \|u - w \|_2^2   \\ [-1.5em]
  &   &  \ \ \ \ \ \ \ \ 
         \ \ \ \ \ \ \ \ + \max_{  s_2 \in \rb^n} -  \frac{1}{2} \| s_2 \|_2^2 
+     s_2^\top \big( \textstyle \frac{u-w}{2} + w_1 - w_2\big)
  \\
 & = & 
  \min_{\substack{w_1 \in \mathcal{W}^{\mathcal{A}_1} \\ w_2 \in \mathcal{W}^{\mathcal{A}_2}}}   
  \bigg(- w_1^\top ( u - w )  + f_1(w_1) + f_2(w_2)
  - \frac{1}{8} \|u - w \|_2^2  \\ [-1.5em]
  &   &  \ \ \ \ \ \ \ \ \ \ \
         \ \ \ \ \ \ \ \ \ \ \ + \frac{1}{2} \| \textstyle \frac{u-w}{2} + w_1 - w_2 \|_2^2 \bigg) \\
& = & 
  \min_{\substack{w_1 \in \mathcal{W}^{\mathcal{A}_1} \\ w_2 \in \mathcal{W}^{\mathcal{A}_2}}}   
  \bigg(- w_1^\top ( u - w )  + f_1(w_1) + f_2(w_2)
  + \frac{1}{2} \|   w_1 - w_2 \|_2^2 \\ [-1.5em]
  &   &  \ \ \ \ \ \ \ \ \ \ \
         \ \ \ \ \ \ \ \ \ \ \ + \frac{1}{2} ( u - w) ^\top ( w_1 - w_2) \bigg)   \\
& = & 
  \min_{\substack{w_1 \in \mathcal{W}^{\mathcal{A}_1} \\ w_2 \in \mathcal{W}^{\mathcal{A}_2}}}   
  \bigg( f_1(w_1) + f_2(w_2) - \frac{1}{2} ( u - w) ^\top ( w_1 + w_2) \\ [-1.5em]
  &   &  \ \ \ \ \ \ \ \ \ \ \
         \ \ \ \ \ \ \ \ \ \ \ + \frac{1}{2} \| w_1 - w_2\|_2^2\bigg),
 \EEAS
 Thus $s_1$ and $s_2$ are dual to certain vectors $w_1$ and $w_2$, which minimize a decoupled formulation in $f_1$ and $f_2$. To check optimality, like in the single function case, it decouples over the constant sets of $w_1$ and $w_2$, which is exactly what step (5) is performing. 
 
 If the check is satisfied, it means that $w_1$ and $w_2$ are in fact optimal for the problem above without the restriction in compatibilities, which implies that they are the Dykstra solutions for the TV problem.
 
 If the check is not satisfied, then the same reasoning as for the one function case, leads directions of descent for the new primal problem above. Hence it decreases; since its value is equal to $-\frac{1}{8} \|s_1 - s_2\|_2^2$, the value of $\|s_1 - s_2\|_2^2$ must increase, hence the second statement.
 
\begin{figure*}
\begin{center}
\begin{tabular}{ccc}
%\hspace*{-0.5cm}
\includegraphics[width=0.3\textwidth]{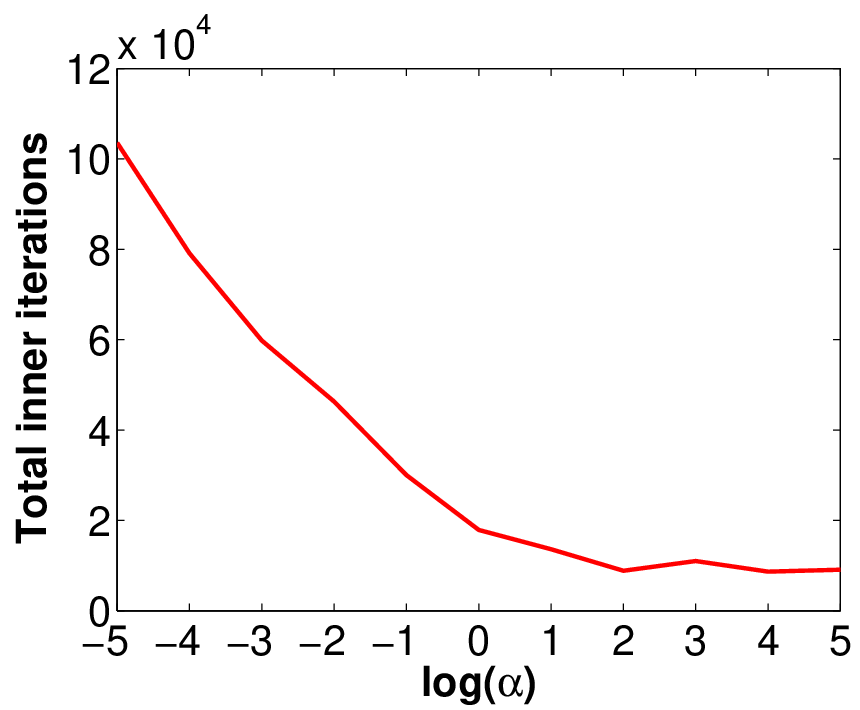} &
%\hspace*{-0.55cm}
\includegraphics[width=0.3\textwidth]{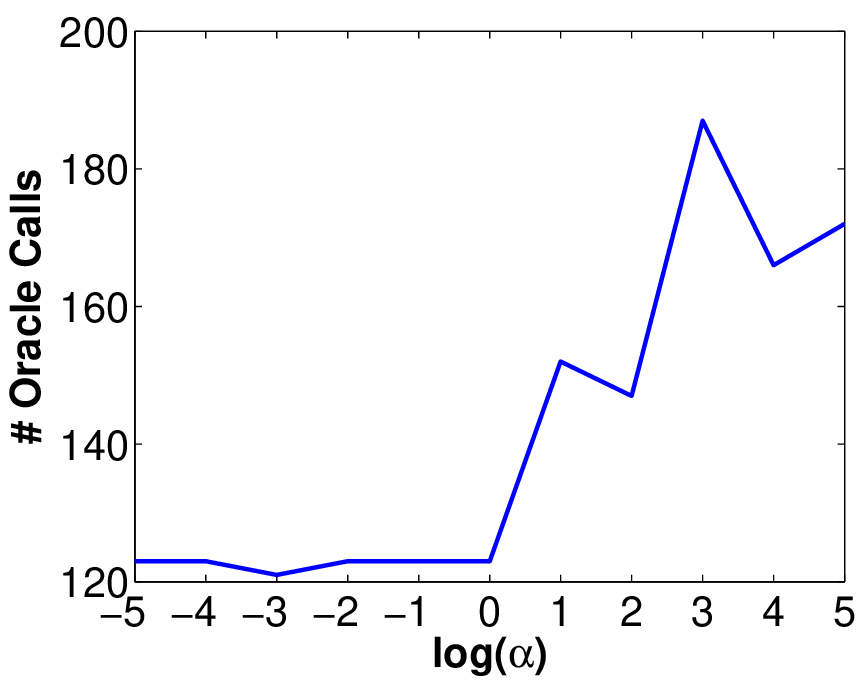} &
%\hspace*{-0.75cm}
\includegraphics[width=0.3\textwidth]{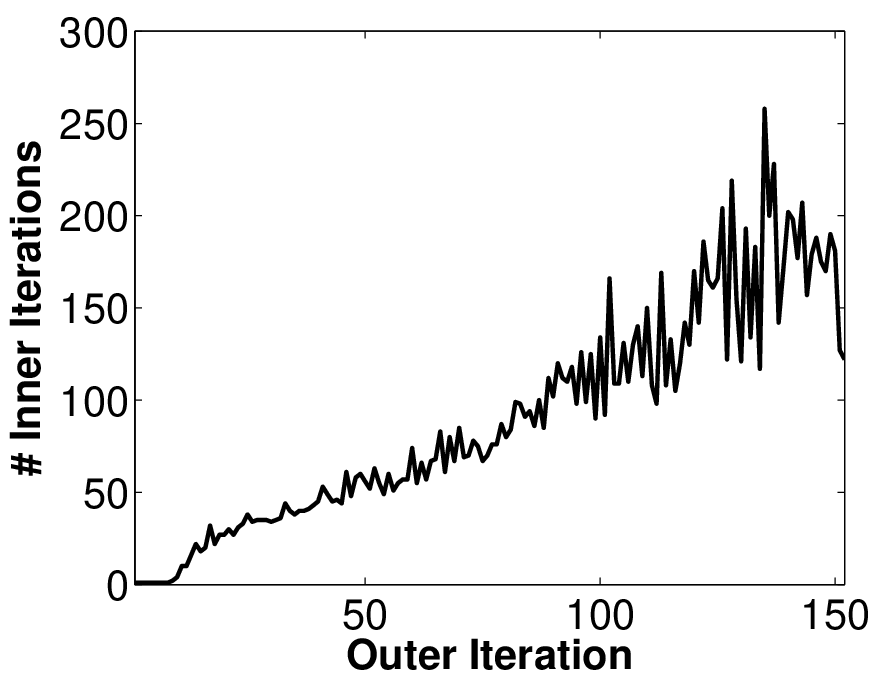}\\
(a) & (b) & (c)\\
\end{tabular}
\end{center}
\caption{(a) Total number of inner iterations for varying $\alpha$. (b) Total number of outer iterations for varying $\alpha$. and (c) Number of inner iterations per each outer iteration for the $\alpha = 10^1$.}
\label{fig:alpha}
\end{figure*}

\section{Decoupled Problems.}
\label{app:decouple}

Given that we deal with polytopes, knowing $w$ implies that we know the faces on which we have to looked for. It turns outs that  for base polytopes, these faces are  products of base polytopes for modified functions (a similar fact holds for  their outer approximations).  

Given the ordered partition $\mathcal{A}'$ defined by the level sets of $w$ (which have to be finer than $\mathcal{A}_1$ and $\mathcal{A}_2$), we know that we may restrict $\widehat{B}^{\mathcal{A}_j}(F_j)$ to elements $s$ such that $s(B) = F(B)$ for all sup-level sets $B$  of $w$ (which have to be unions of contiguous elements of $A_j$); see an illustration below.

\begin{center}
\includegraphics[scale=.375]{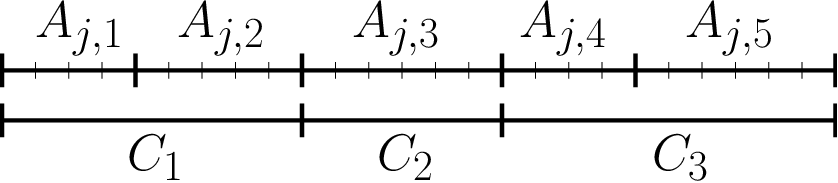}  
\end{center}

More precisely, if $C_1,\dots,C_{m'}$ are constant sets of $w$ ordered with decreasing values. Then, we may search for $s_j$ independently for each subvector $(s_j)_{C_k} \in \rb^{C_k}$, $ k \in \{1,\dots,m'\}$ and with the constraint that
$$(s_j)_{C_k} \in \hat{B}^{\mathcal{A}_j \cap C_k}\big[ ( F_ j) _{ C_k |  C_1 \cup \cdots \cup C_{k-1} }\big],$$
where  $\mathcal{A}_j \cap C_k$ is the ordered partition obtained from $\mathcal{A}_j$ once restricted onto $C_k$ and the submodular function is the so-called contraction of $F$ on $C_k$ given $C_1 \cup \cdots \cup C_{k-1}$, defined as $ S \mapsto  F_j( S \cup C_1 \cup \cdots \cup C_{k-1} ) - 
F(C_1 \cup \cdots \cup C_{k-1}) $. Thus this corresponds to solving $m$ different smaller subproblems.

\section{Approximate Dykstra Steps}
\label{app:alpha}

The Dykstra step, i.e., step (4) of the algorithm proposed in \mysec{Dykstra} is not finitely convergent. Therefore, it needs to be solved approximately.  For this purpose, we introduce a parameter $\alpha$ to approximately solve the Dykstra step such that $\| s_1 + s_2 - u + w \|_1 \leq \alpha (\epsilon_1 + \epsilon_2)$.  Let $\epsilon$ be defined as $\alpha ( \epsilon_1 + \epsilon_2)$. This shows that the $s_1$ and $s_2$ are $\epsilon$-accurate. Therefore, $\alpha$ must be chosen in such a way that we avoid cycling in our algorithm. However, another alternative is to warm start the Dykstra step with $w_1$ and $w_2$ of the previous iteration. This ensures we don't go back to the same $w_1$ and $w_2$, which we have already encountered and avoid cycling. \myfig{alpha} shows the performance of our algorithm for a simple problem of $100 \times 100$ 2D-grid with 4-neighborhood and uniform weights on the edges with varying $\alpha$. \myfig{alpha}-(a) shows the total number of inner iterations required to solve the TV problem. \myfig{alpha}-(b) gives the total number of SFM oracle calls required to solve the TV problem. In \myfig{alpha}-(c), we show the number of inner iterations in every outer iteration for the best $\alpha$ we have encountered.

\bibliography{super_tree}

\end{document}